\documentclass[a4paper,11pt,twoside,reqno]{amsart}
\usepackage{amssymb,amsmath,amsfonts,amsthm}
\usepackage{color}
\usepackage[dvipsnames]{xcolor}
\usepackage{graphicx}
\usepackage{placeins}
\usepackage{tikz}
\usepackage{pgfplots}
\pgfplotsset{compat=1.18}
\usepackage[margin=1in]{geometry}
\usepackage{hyperref}
\newtheorem{theorem}{Theorem}[section]
\newtheorem{lemma}[theorem]{Lemma}
\newtheorem{proposition}[theorem]{Proposition}
\newtheorem{corollary}[theorem]{Corollary}
\newtheorem{remark}[theorem]{Remark}

\numberwithin{equation}{section}

\renewcommand{\epsilon}{\varepsilon}
\newcommand{\vol}{\mathrm{vol}} 

\newcommand{\RR}{\mathbb{R}}
\newcommand*{\CC}{\mathbb{C}}

\newcommand*{\NN}{\mathbb{N}}

\newcommand{\lquot}[2]{\raisebox{-0.5ex}{$#2$}\backslash\!\raisebox{0.5ex}{$#1$}}
\newcommand{\rquot}[2]{\raisebox{0.5ex}{$#1$}\!/\!\raisebox{-0.5ex}{$#2$}}

\title{On the spectrum of the magnetic Dirac operator}

\author{Volker Branding}
\address{University of Rostock, Institute of Mathematics\\
Ulmenstraße 69, 18057 Rostock, Germany}
\email{volker.branding@uni-rostock.de}

\author{Nicolas Ginoux}
\address{Universit\'{e} de Lorraine, CNRS, IECL, F-57000 Metz, France}
\email{nicolas.ginoux@univ-lorraine.fr}

\author{Georges Habib}
\address{Lebanese University, Faculty of Sciences {\bf II}, Department of Mathematics, 
P.O. Box 90656 Fanar-Matn, Lebanon and Universit\'e de Lorraine, CNRS, IECL, F-54506, Nancy, France}
\email{ghabib@ul.edu.lb}


\date{\today}
\begin{document}

\begin{abstract}
The magnetic Dirac operator describes the relativistic motion of a charged particle in a magnetic field.
Although this operator got a lot of attention in physics many of its fundamental mathematical properties remain unexplored and
this article is a first step towards filling this gap. 
To this end we provide a number of eigenvalue estimates for the magnetic Dirac operator
on closed Riemannian manifolds and explicitly compute its spectrum for specific choices of the magnetic field on the flat torus and on the three-dimensional round sphere.
\end{abstract}

\maketitle

\section{Introduction}\label{s:intro}

The Dirac operator is a first-order differential operator which was originally introduced by physicist Paul Dirac to describe electrons, which are spin $\frac{1}{2}$-particles, in spacetime.
In the presence of an electromagnetic field, a magnetic potential has to be added to the Dirac operator, turning it into the so-called \emph{magnetic Dirac operator}.
Over the last few years, a lot of studies have been devoted to the magnetic Dirac operator, mainly from the perspectives of physics and analysis, see e.g. \cite[Sec. 1]{CharalambousGrosse:23} as well as \cite[Sec. 1]{FrankLoss:2022} and references therein for an overview.
In particular, there has been strong interest in the so-called zero modes, which are the elements of the kernel of the magnetic Dirac operator, see \cite{FrankLoss:2021,FrankLoss:2022,Reuss:25}, where \cite{Reuss:25} extends the results \cite{FrankLoss:2021,FrankLoss:2022} from flat space to arbitrary closed Riemannian spin manifolds.\\

In this article, we make a decisive contribution to the spectral theory of the magnetic Dirac operator on closed Riemannian spin manifolds.
We mainly focus on the interactions between the spectrum of that operator and the geometry of the underlying manifold.
We derive general geometric estimates for the smallest eigenvalues of the magnetic Dirac operator and compute its whole spectrum in two cases.
These estimates generalise Friedrich's and B\"ar's resp. Hijazi's ones \cite{Baer:92,Friedrich:80,Hijazi:86} and their limiting cases restrict a lot the geometry of the manifold when the magnetic field is nowhere zero. \color{black}
We also keep track of the spectral shift occurring when going from the Dirac to the magnetic Dirac operator.
Surprisingly enough, introducing a magnetic field, even a small one, does not necessarily make the first eigenvalue decrease, as one would expect from the so-called \emph{diamagnetic inequality} on the connection level \cite[Sec. 2.4]{CharalambousGrosse:23}.\\

The article is structured as follows.
After a review of the magnetic Dirac operator in Section \ref{sec:review}, we state and prove two general but fundamental estimates for its smallest eigenvalue in Section \ref{sec:estimates}, see Theorems \ref{thm:magneticfrie}
 and \ref{t:Hijazimagneticineq}. 
The equality case of those estimates is carefully studied in the particular case when the magnetic field is nowhere zero. 
 As a result, it turns out that the geometry of the manifold restricts to a Sasaki structure and the magnetic vector field is the corresponding Reeb vector field. 
 This motivates the study in Section \ref{s:KillingmagneticVF} when a magnetic field is a Killing vector field of constant norm. \color{black}
On the way, we notice that, thanks to the conformal covariance of the Dirac operator, the proof of the main inequality in \cite[Theorem 1]{Reuss:25} can be considerably simplified, see Remark \ref{rem:proofReuss}.
On surfaces, a fine estimate taking into account the particular structure of the vanishing set of an eigenvector for the magnetic Dirac operator can be established, see Theorem \ref{t:estimnodalset}.
The difference between the respective first eigenvalues of the standard and magnetic Dirac operator is studied in Section \ref{s:diamagneticinequality}, where we focus on the diamagnetic inequality comparing both eigenvalues.
The main tool for testing that inequality is set up in Proposition \ref{prop:diamine} and relies on the use of the min-max principle applied to an eigenspinor of the standard Dirac operator as a test spinor. 
We examine this inequality on Einstein-Sasakian manifolds. \color{black} 
Section \ref{s:computmagneticspectra} is devoted to the computation of the complete spectrum of the magnetic Dirac operator on both the $3$-dimensional sphere with standard Sasaki structure and on the flat $n$-dimensional torus, see Theorems \ref{t:spectrumDetaS3} and \ref{t:specflattorusspincAparallel}. 
Thanks to this computation, we show that the diamagnetic inequality can never hold on the sphere and, depending on the choice of the magnetic field, it may or may not hold on the flat torus.\color{black}
The particular case where the magnetic field is Killing of constant length is the topic of Section \ref{s:KillingmagneticVF}. 
In this particular setting, the manifold is locally submerged into a base manifold whose fibers are just the integral curves of the magnetic vector field. 
A natural question which arises in this setup is how the spectrum of the magnetic Dirac operator can be expressed in terms of the geometry of those submersions. 
For this, we \color{black} consider the magnetic Dirac eigenvalues that have so-called basic associated eigenspaces and show that they can be bounded in a much finer way than in Theorem \ref{thm:magneticfrie}, see Theorem \ref{thm:estimatebasic}.\\

We underline that the presence of the magnetic field does not reduce to a technical difficulty when comparing with the standard Dirac operator.
It makes its spectral behaviour pointwise very different, as is exemplified in the diamagnetic inequality mentioned above.\\

{\bf Acknowledgment:} This work was carried out during the visit of the third named author to University of Lorraine as a visiting Professor in November 2024 and July 2025. 
He would like to thank IECL and IHES for the hospitality. 
The second and the third named authors are also supported by the International Emerging Action ``HOPF" (Higher Order boundary value Problems for differential Forms) of the French CNRS, which they would like to thank.

\section{Review on the magnetic Dirac operator}\label{sec:review}
In this section, we review some basic facts on the magnetic Dirac operator defined on a spin manifold. 
For more details, we refer to \cite{CharalambousGrosse:23,Savale:17,Savale:18}.\\

Let $(M,g)$ be a Riemannian spin manifold of dimension $n$ and let $\eta$ be a real $1$-form on $M$. 
We denote by $\Sigma M$ its complex spinor bundle. 
Recall that this bundle is equipped with a Hermitian product $\langle\cdot\,,\cdot\rangle$ and a metric connection $\nabla$ coming from the Levi-Civita connection on $(M,g)$ such that the Clifford action of any vector field on sections of $\Sigma M$ is skew-Hermitian and parallel. 
Sections of $\Sigma M$ are called spinors. 
The Clifford action of vector fields on spinors satisfies the so-called Clifford relations, which read $X\cdot Y\cdot\varphi+Y\cdot X\cdot\varphi=-2g(X,Y)\varphi$, for all vector fields $X,Y$ and any spinor $\varphi$ on $M$.
Recall also that the Dirac operator is the first order differential operator defined as $\displaystyle D:=\sum_{k=1}^n e_k\cdot\nabla_{e_k}$, where $\{e_k\}_{k=1,\ldots,n}$ is any local orthonormal frame of $TM$. 
The magnetic Dirac operator $D^\eta$ acts on the spinor bundle $\Sigma M$ by the following
$$D^\eta:=D+i\eta\cdot.$$
It is not difficult to check that the magnetic Dirac operator is an elliptic and essentially self-adjoint operator when $M$ is complete. 
Thus, when $M$ is closed, the magnetic Dirac operator has a discrete spectrum made of real eigenvalues of finite multiplicities, which we denote by $(\lambda_k^\eta)_k$, and the corresponding eigenspaces consist of smooth sections only.
Let us define a new metric connection on $\Sigma M$ by the following $\nabla^\eta_X:=\nabla_X+i\eta(X)$. 
It is easy to check that this connection is compatible with the Clifford multiplication, that is $\nabla_X^\eta(Y\cdot)=(\nabla_X Y)\cdot+Y\cdot\nabla_X^\eta$ for any vector fields $X,Y\in TM$ and the magnetic Dirac operator can be written in terms of $\nabla^\eta$ as $D^\eta=\sum_{k=1}^n e_k\cdot\nabla^\eta_{e_k}$. 
When we identify the Clifford module $\mathbb{C}\mathrm{l}(TM)$ with the exterior algebra, we can write $D^\eta=d^\eta+\delta^\eta$, with $d^\eta=d+i\eta\wedge$ and $\delta^\eta=\delta-i\eta^\sharp\lrcorner$.
In this case, $(D^\eta)^2$ is in general not equal to the magnetic Hodge Laplacian $\Delta^\eta:=d^\eta\delta^\eta+\delta^\eta d^\eta$, since $(d^\eta)^2=id\eta\wedge$ which does not vanish in general \cite{EGHP:23}.\\   

On the other hand, the magnetic Dirac operator can be seen as the Dirac operator of a spin$^c$ structure whose auxiliary line bundle is trivial and carrying a connection given by $2i\eta$. 
In particular, when $M$ is simply connected and $d\eta=0$, the spin$^c$ connection can be identified with the spin connection \cite[Lem.2.1]{Moroianu:97}, that is, we have $\nabla^\eta=e^{-if}\nabla(e^{if})$ where $\eta=df$. 
Now, one can easily check that the curvature $R^\eta$ associated with the connection $\nabla^\eta$ is given by 
\begin{equation}\label{eq:curvatureeta}
R^\eta(X,Y)=R(X,Y)+i(d\eta)(X,Y)
\end{equation}
for any vector fields $X, Y\in TM$ where $R$ is the curvature associated with $\nabla$. 
Based on this identity, the magnetic Schr\"odinger-Lichnerowicz formula for $(D^\eta)^2$ can be stated as follows:

\begin{proposition}\label{p:SchroedingerLichnerowiczmagnetic} 
We have
\begin{align}
\label{eq:SchroedingerLichnerowiczmagnetic}
    (D^\eta)^2=(\nabla^\eta)^*\nabla^\eta+\frac{1}{4}S+id\eta\cdot,
\end{align}
where $S$ is the scalar curvature of the manifold $(M,g)$. 
Also, we have
\begin{eqnarray}\label{eq:relationdiracmagnwithdira}
(D^\eta)^2&=&D^2+id\eta\cdot+i(\delta\eta)\cdot-2i\nabla_{\eta}+|\eta|^2.
\end{eqnarray}
\end{proposition}

\begin{proof}
Using the expression of the magnetic Dirac operator, we choose an orthonormal frame $\{e_k\}_k$ which is parallel at some point $x\in M$ and compute, at $x$,
\begin{eqnarray*}
(D^\eta)^2&=&\sum_{k,l=1}^n e_k\cdot\nabla^\eta_{e_k} (e_l\cdot \nabla^\eta_{e_l}) \\
&=&\sum_{k,l=1}^n e_k\cdot e_l\cdot\nabla^\eta_{e_k} \nabla^\eta_{e_l}\\
&=&-\sum_{k=1}^n \nabla^\eta_{e_k} \nabla^\eta_{e_k}+\frac{1}{2}\sum_{k,l=1}^ne_k\cdot e_l\cdot R^\eta(e_k,e_l)\\
&\stackrel{\eqref{eq:curvatureeta}}{=}&-\sum_{k=1}^n \nabla^\eta_{e_k} \nabla^\eta_{e_k}+\frac{1}{2}\sum_{k,l=1}^ne_k\cdot e_l\cdot R(e_k,e_l)+\frac{i}{2}\sum_{k,l=1}^n(d\eta)(e_k,e_l)e_k\cdot e_l\cdot\\
&=&(\nabla^\eta)^*\nabla^\eta+\frac{1}{4}S+id\eta\cdot.
\end{eqnarray*}
Identity \eqref{eq:relationdiracmagnwithdira} comes from expanding $(D^\eta)^2=(D+i\eta\cdot)^2$ along with $D(\eta\cdot)=(d\eta+\delta\eta)\cdot-2\nabla_\eta-\eta\cdot D$. 
This ends the proof.
\end{proof}

In the following, we show the gauge invariance of the magnetic Dirac operator
\begin{proposition} \label{prop:gaugeonva} Let $(M^n,g)$ be a Riemannian spin manifold and let $\eta$ be a differential form on $M$. 
For any $\eta_\tau:=\frac{d\tau}{i\tau}$ with $\tau\in C^\infty(M,\mathbb{U}_1)$, the magnetic Dirac operators $D^\eta$ and $D^{\eta+\eta_\tau}$ are unitarily equivalent, meaning that
$$D^{\eta+\eta_\tau}=\overline{\tau}D^\eta\tau.$$
In particular, if $M$ is closed, the operators $D^\eta$ and $D^{\eta+df}$ are isospectral for any smooth real-valued function $f$.
\end{proposition}
\begin{proof} For any spinor field $\psi$, we write
\begin{eqnarray*}
D^\eta(\tau\psi)&=& D(\tau\psi)+i\tau\eta\cdot\psi\\
&=&\tau D^\eta\psi+d\tau\cdot\psi\\
&=&\tau\left(D^\eta\psi+i\frac{d\tau}{i\tau}\cdot\psi\right)\\
&=&\tau D^{\eta+\eta_\tau}\psi.
\end{eqnarray*}
For the last part, we just take $\tau=e^{if}$. 
This finishes the proof.
\end{proof}

\begin{proposition}\label{p:symmetryspecDiracevendim}
If the dimension $n$ of a Riemannian spin manifold $(M^n,g)$ is even, then for any $\eta\in\Omega^1(M)$, the magnetic Dirac operator $D^\eta$ anti-commutes with the Clifford action of the volume form of $M$.
As a consequence, if $M$ is closed, then $D^\eta$ has symmetric spectrum.
\end{proposition}

\begin{proof}
It suffices to notice that the Riemannian volume form provided by the orientation of $M$ not only anti-commutes with the Dirac operator of $(M^n,g)$ but also with the Clifford action of $\eta$.
Therefore, it anti-commutes with $D^\eta$.
\end{proof}

\section{Eigenvalue estimates for the magnetic Dirac operator}\label{sec:estimates}
In this section, we  establish some eigenvalue estimates {\it \`a la Friedrich} \cite{Friedrich:80} and {\it \`a la  B\"ar} \cite{Baer:92} resp. {\it \`a la Hijazi} \cite{Hijazi:86} for the magnetic Dirac operator and discuss their limiting cases. 
It turns out that, when equality holds in these estimates, the magnetic field $\eta$ gives rise to a particular geometry on the manifold which does not necessarily reduce to the case without magnetic field, see the discussions in Theorems \ref{thm:magneticfrie} and \ref{t:Hijazimagneticineq}. \\

Before stating these results, let us recall some properties of the spinor bundle on a Sasakian manifold (see \cite{FriedrichKim:00} for more details).  
Given a Sasakian spin manifold $(M,g,\eta)$ of dimension $n=2m+1$ with Reeb vector field $\eta$, the spinor bundle of $M$ decomposes under the action of the transversal K\"ahler form  $\Omega=\frac{1}{2}d\eta$ into
\begin{equation}\label{eq:decomspinor}
\displaystyle\Sigma M=\bigoplus_{r=0}^m\Sigma_r M,
\end{equation}
where $\Sigma_r M$ is the eigenbundle associated with the eigenvalue $i(2r-m)$ of $\Omega$.  
The Clifford action of $\eta$ on $\Sigma_r M$ is given by $\eta\cdot_{|_{\Sigma_r M}}=i(-1)^{r+m} {\rm Id}_{|_{\Sigma_r M}}$. 
A Sasakian manifold $(M^{2m+1},g,\eta)$ is called $\eta$-Einstein if the Ricci curvature satisfies 
$${\rm Ric}=\alpha g+\beta\eta\otimes\eta,$$
for some $\alpha,\beta\in C^\infty(M)$. 
It was shown in \cite{Okumura:62} that, if $m>1$, the functions $\alpha$ and $\beta$ are constants satisfying $\alpha+\beta=2m$ and, in this case, the scalar curvature is equal to $2m(\alpha+1)$.
On a Sasakian spin manifold the notion of {\it Sasakian quasi-Killing spinor of type $(a,b)$} as being a solution of the differential equation
$$\nabla_X\psi=aX\cdot\psi+b\eta(X)\eta\cdot\psi,$$
for real numbers $a$ and $b$ was defined in \cite{FriedrichKim:00}.
Moreover, in \cite[Lem. 6.5]{FriedrichKim:00} it is proven that the existence of a nonzero Sasakian quasi-Killing spinor of type $(\pm\frac{1}{2},b)$ with $b\neq 0$ implies that the manifold is $\eta$-Einstein of constant $\alpha=2m\pm 4b$.
The case when $b=0$ corresponds to real Killing spinors. Furthermore, \cite[Thm. 6.3]{FriedrichKim:00} and \cite[Thm. 8.4] {FriedrichKim:00} show that any simply connected $\eta$-Einstein Sasakian manifold admits a Sasakian quasi-Killing spinor $\psi_m\in \Sigma_m M$ of type $(\frac{-1}{2},b)$. 
When $m$ is odd, it also has a Sasakian quasi-Killing spinor $\psi_0\in\Sigma_0 M$ of the same type. Note that such quasi-Killing spinors will play a crucial role in the characterization of the equality case of the magnetic Friedrich inequality. With the help of the magnetic Schr\"odinger-Lichnerowicz-formula from Proposition \ref{p:SchroedingerLichnerowiczmagnetic}, we state the magnetic Friedrich inequality (see \cite{HerzlichMoroianu:99} and \cite[Thm. 2.3]{Roos:16} for the corresponding inequalities on spin$^c$ manifolds). 

\begin{theorem} \label{thm:magneticfrie} (Magnetic Friedrich Inequality) Let $(M^n,g)$ be a closed Riemannian spin manifold and let $\eta\in \Omega^1(M)$. 
For any $t\in \mathbb{R}_+$, any eigenvalue $\lambda^{t\eta}$ of the magnetic Dirac operator $D^{t\eta}$ satisfies 
$$(\lambda^{t\eta})^2\geq \frac{n}{4(n-1)}\mathop{\rm inf}\limits_{M^n}\big(S-4t\lfloor\frac{n}{2}\rfloor^\frac{1}{2}|d\eta|\big).$$
If equality occurs, then either we are in the equality case of Friedrich's inequality (i.e. $\eta=0$) or, up to rescaling the metric, $\lambda^{t\eta}=\pm \frac{n}{2}$ and the universal cover $\widetilde M$ of $M$ is a non-Einstein Sasakian manifold. 
If furthermore $\eta^{\sharp_g}$ is a nontrivial geodesic vector field of constant norm, which thus can be assumed to be equal to $1$ up to rescaling $\eta$, and $M$ is simply-connected,
equality is realized if and only if $M$ is an $\eta$-Einstein non-Einstein Sasakian manifold of constant scalar curvature $(n-1)(n+4t)$.  
\end{theorem} 

\begin{proof}
To prove the estimate, we proceed as in the case without magnetic field. 
For this, we define a magnetic twistor operator by
\begin{align*}
P^\eta_X\psi:=\nabla^\eta_X\psi+\frac{1}{n}X\cdot D^\eta\psi    
\end{align*}
for all \(X\in TM\).
This operator satisfies the pointwise equality
\begin{align*}
|\nabla^\eta\psi|^2=|P^\eta\psi|^2+\frac{1}{n}|D^\eta\psi|^2.    
\end{align*}

Together with the magnetic Schrödinger-Lichnerowicz formula \eqref{eq:SchroedingerLichnerowiczmagnetic}
we  arrive at
\[\int_M\left(|D^{t\eta}\psi|^2-\frac{n}{4(n-1)}\left(S|\psi|^2
     +4t\langle i d\eta\cdot\psi,\psi\rangle
    \right)\right)d\mu_g=\frac{n}{n-1}\int_M|P^{t\eta}\psi|^2\,d\mu_g\geq 0.
\]

Recall that in \cite[Lemma 3.3]{HerzlichMoroianu:99} the following estimate 
\begin{align*}
    \langle i \Omega\cdot\psi,\psi\rangle\geq -\lfloor\frac{n}{2}\rfloor^\frac{1}{2}
    |\Omega||\psi|^2
\end{align*}
was established for any \(\psi\in\Gamma(\Sigma M)\) and any differential two-form \(\Omega\). 
Here the norm of $\Omega$ is considered as the norm of a differential two-form, that is $\displaystyle|\Omega|^2=\sum_{k<l}\Omega(e_k,e_l)^2$. 
Equality is attained when $\psi\neq0$ if and only if either $\Omega$ vanishes or has maximal rank equal to $n$ if $n$ is even or to $n-1$ if $n$ is odd.  
Hence, if \(\psi\) is an eigenspinor for the magnetic Dirac operator associated with the eigenvalue $\lambda^{t\eta}$, we find
\[\int_M\left((\lambda^{t\eta})^2-\frac{n}{4(n-1)}\left(S
     -4t\lfloor\frac{n}{2}\rfloor^\frac{1}{2}
    |d\eta|
    \right)\right)|\psi|^2d\mu_g\geq 0
    \]
leading to the magnetic Friedrich inequality. Assume now that equality is attained, then we have equality in all above inequalities. Therefore, the spinor field $\psi$ satisfies 
\begin{equation}\label{eq:eqkillin}
\nabla_X^{t\eta}\psi=-\frac{\lambda^{t\eta}}{n} X\cdot\psi\quad\text{and}\quad d\eta\cdot\psi=i\lfloor\frac{n}{2}\rfloor^\frac{1}{2}
    |d\eta|\psi.
    \end{equation}
Here, either $d\eta=0$ or $d\eta$ has maximal rank as stated before. 
Recall that $\psi$ corresponds to a Killing or parallel spinor for the spin$^c$ structure with trivial auxiliary bundle of curvature $\Omega=2itd\eta$.  Thus as mentioned in Section \ref{sec:review}, when $d\eta=0$, the spin$^c$ connection on the universal cover of $M$ corresponds to the spin connection and thus we are in the equality case of the usual Friedrich inequality. 
Hence, we are left with the case when $d\eta$ has maximal rank. In the following, we distinguish two cases: the case when $\lambda^{t\eta}=0$, meaning that $\psi$ is parallel for the magnetic connection $\nabla^{t\eta}$, or $\lambda^{t\eta}\neq 0$ meaning that $\psi$ is Killing. \\

Let us first discuss the case where $\psi$ is $\nabla^{t\eta}$-parallel. 
When $\nabla^{t\eta}\psi=0$, the magnetic Schr\"odinger-Lichnerowicz formula gives that $itd\eta\cdot\psi=-\frac{1}{4}S\psi$. 
Therefore from \eqref{eq:eqkillin}, we deduce that $S=4t\lfloor\frac{n}{2}\rfloor^\frac{1}{2}|d\eta|\geq 0$.  Hence, we get that $ S=2\lfloor\frac{n}{2}\rfloor^\frac{1}{2}|\Omega|$. Therefore by \cite[Prop. 3.3]{EHGN:17}, the universal cover of $M$ is isometric to either a spin manifold with parallel spinors, a K\"ahler-Einstein manifold of nonnegative scalar curvature or the Riemannian product of a K\"ahler-Einstein manifold of nonnegative scalar curvature with $\mathbb{R}$. Now, the last two cases cannot occur since K\"ahler-Einstein manifolds cannot have Ricci form equal to $2itd\eta$. Indeed, the Ricci form of a K\"ahler-Einstein manifold is equal, up to some constant, to the K\"ahler form of that manifold which cannot be an exact form. Hence the universal cover should be Ricci-flat and thus, $S=0$. 
Therefore $\Omega$ must vanish, i.e.
$d\eta=0$, which contradicts the fact that $d\eta$ is of maximal rank.
This shows that the case of $\psi$ being $\nabla^{t\eta}$-parallel cannot occur.\\

Let us now discuss the case when $\lambda^{t\eta}\neq 0$.
Hence, by rescaling the metric into the form $\overline{g}=c^2g$ where $c=2\frac{\lambda^{t\eta}}{n}$, we obtain $\lambda^{t\eta}\in\{\pm\frac{n}{2}\}$.
Namely taking a metric of the form $\overline{g}=c^2 g$ for some constant $c$ to be determined, the first equation in \eqref{eq:eqkillin} becomes with respect to the metric $\overline{g}$ equal to 
$$\overline\nabla_X\overline{\psi}=\overline{\nabla_X\psi}=-\frac{\lambda^{t\eta}}{n}\overline{X}\overline\cdot\overline{\psi}-it\eta(X)\overline\psi=-\frac{\lambda^{t\eta}}{cn}X\overline\cdot\overline{\psi}-it\eta(X)\overline\psi.$$
By choosing $c$ such that $-\frac{\lambda^{t\eta}}{cn}=\pm \frac{1}{2}$,
the equation then reduces to 
\begin{equation}\label{eq:nablapsi}
\overline\nabla_X\overline{\psi}=\pm \frac{1}{2}X\overline\cdot~\overline{\psi}-it\eta(X)\overline\psi.
\end{equation}
Now a simple computation gives that
$|d\eta|_g=c^2 |d\eta|_{\overline{g}}$. 
Therefore the magnetic Friedrich inequality -- and thus its equality case as well -- remain unchanged.
Hence, the second equation in \eqref{eq:eqkillin} becomes 
\begin{equation}\label{eq:detapsi}
d\eta\overline{\cdot}~\overline\psi=i\lfloor\frac{n}{2}\rfloor^\frac{1}{2}|d\eta|_{\overline g}\overline{\psi}
\end{equation}
where the action of $d\eta$ is given by 
$d\eta\overline{\cdot}=\sum_{i,j} d\eta(\overline{e_i},\overline{e_j})\overline{e}_i~{\overline \cdot}~
\overline{e}_j~{\overline \cdot}$. 
The equality case together with $\lambda^{t\eta}\in\{\pm\frac{n}{2}\}$ yields $\overline{S}-4t\lfloor\frac{n}{2}\rfloor^\frac{1}{2}|d\eta|_{\overline g}=n(n-1)$. 
From \cite[Thm. 4.1]{Moroianu:97}, the universal cover $(\widetilde M,\overline{g})$ of $M$ must be a Sasakian manifold. 
Notice here that $\widetilde M$ cannot be Einstein, since otherwise this would imply $d\eta=0$ \cite[Prop. 4.2]{Moroianu:97} which would again contradict the maximality of the rank.  
This proves the first part of the theorem.
To show the second part, where we assume from now on that $\eta$ has constant length, we first notice that $\eta^{\sharp_{\overline{g}}}=\frac{1}{c^2}\eta^{\sharp_{g}}$, and therefore, $\eta^{\sharp_{\overline{g}}}$ will stay geodesic of constant norm with respect to the metric $\overline{g}$. 
Up to replacing $\eta$ by $\frac{\eta}{|\eta|}$ and $t$ by $|\eta|t$, we may assume that $\eta$ is of unit length.
To simplify the notations, we will denote $(M,\overline{g})$ by $(M,g)$. Taking the derivative of \eqref{eq:nablapsi}, the curvature of the spinor $\psi$ is equal to  $R(X,Y)\psi=\frac{1}{4}(Y\cdot X-X\cdot Y)\cdot\psi-itd\eta(X,Y)\psi$. 
Hence by tracing and using the Ricci identity formula $-\frac{1}{2}{\rm Ric}(X)\cdot\psi=\sum_i e_i\cdot R(X,e_i)\psi$ (see e.g. \cite[Lemma 1.2.4]{Ginoux:2009}), we deduce that
\begin{equation}\label{eq:riccidentity}
{\rm Ric}(X)\cdot\psi=(n-1)X\cdot\psi+2it(X\lrcorner d\eta)\cdot\psi,
\end{equation} 
for all $X\in TM$. 
As $\widetilde{M}$ is Sasakian, we denote by $\widetilde\zeta$ the corresponding Reeb vector field. 
Since $\widetilde{\rm Ric}(\widetilde\zeta)=(n-1)\widetilde\zeta$, we get from \eqref{eq:riccidentity} that $\widetilde\zeta\lrcorner d\widetilde\eta=0$, where $\widetilde\eta$ is the pull-back of $\eta$ to $\widetilde M$. 
On the other hand, we have from the fact that $\widetilde\eta^{\sharp_g}$ is geodesic of constant norm that $\widetilde\eta^{\sharp_g}\lrcorner d\widetilde\eta=0$. 
Hence, as the rank of $d\eta$ (and, thus, of $d\widetilde\eta$) is maximal equal to $n-1$ because of $n=2m+1$ being odd, we deduce that $\widetilde\zeta$ must be parallel to $\widetilde\eta$. 
Therefore $\widetilde\zeta=\pm\widetilde\eta^{\sharp_g}$.
Up to replacing $\widetilde\zeta$ by $-\widetilde\zeta$, we may assume that $\widetilde\zeta=\widetilde\eta^{\sharp_g}$.
Using Equation \eqref{eq:detapsi} and the fact that $\widetilde M$ is Sasakian of dimension $n=2m+1$, we deduce that 
$$d\widetilde\zeta\cdot\psi=i\lfloor\frac{n}{2}\rfloor^\frac{1}{2}|d\widetilde\zeta|_g\psi=2im\psi,$$
and, thus, $\widetilde\Omega\cdot\psi=im\psi$, where $\widetilde\Omega=\frac{1}{2}d\widetilde\zeta$ is the transversal K\"ahler form of the Sasakian manifold $\widetilde M$. 
Here, we use the fact that the norm of the transversal K\"ahler form is equal to 
$$|\widetilde\Omega|^2=\sum_{k<l}\widetilde\Omega(e_k,e_l)^2=\frac{1}{2}\sum_{k,l}g(Je_k,e_l)^2=m,$$
where $\{e_k\}_{k=1,\ldots,2m}$ is a local orthonormal frame of $\widetilde\zeta^\perp$. 
Hence, according to the decomposition \eqref{eq:decomspinor}, we deduce that $\psi\in \Sigma_m\widetilde M$. 
Therefore, by the identities $(X\lrcorner\widetilde\Omega)\cdot\psi=-iX\cdot\psi$ for all $X\perp\widetilde\zeta$ and $\widetilde\zeta\cdot\psi=i\psi$ (for $\psi\in\Sigma_m \widetilde M$), we obtain $\frac{1}{2}(X\lrcorner d\widetilde\zeta)\cdot\psi=-iX\cdot\psi-\widetilde\zeta(X)\psi$ for all tangent vectors $X$, including $X=\widetilde\zeta$.
Hence \eqref{eq:riccidentity} becomes 
\begin{eqnarray*}
\widetilde{\rm Ric}(X)\cdot\psi&=&(n-1)X\cdot\psi+2it(-2iX\cdot\psi-2\widetilde\zeta(X)\psi)\\
&=& ((n-1+4t)X-4t\widetilde\zeta(X)\widetilde\zeta)\cdot\psi,
\end{eqnarray*}
for any $X\in T\widetilde M$. 
Therefore, the manifold $\widetilde M$ is $\eta$-Einstein with Ricci tensor equal to $\widetilde{\rm Ric}=(n-1+4t)g-4t\widetilde\zeta\otimes\widetilde\zeta$ and, thus,  of constant scalar curvature equal to $S=(n-1)(n+4t)$.
Note that $\widetilde M$ must be compact because of $\widetilde{\rm Ric}\geq(n-1)g$.\\\\
 
To check the converse, we let $(M^{2m+1},g,\eta)$ to be any simply connected closed $\eta$-Einstein Sasakian spin manifold of constant scalar curvature equal to $2m(2m-4b+1)$ with $b\in \mathbb{R}$. 
It is shown in \cite[Thm. 6.3]{FriedrichKim:00} and \cite[Thm. 8.4]{FriedrichKim:00} that the manifold $M$ admits a Sasakian quasi-Killing spinor $\psi_m\in \Sigma_m M$ of type $(\frac{-1}{2},b)$. 
Hence $\psi_m$ is an eigenspinor for the Dirac operator associated with the eigenvalue $\frac{2m+1}{2}-b$ and we have 
$$D^{t\eta}\psi_m=D\psi_m+it\eta\cdot \psi_m=\left(\frac{2m+1}{2}-b-t\right)\psi_m,$$
where we use that $\eta\cdot_{|_{\Sigma_r M}}=i(-1)^{r+m} {\rm Id}_{|_{\Sigma_r M}}$. 
Hence for $b=-t$, we get that $\frac{2m+1}{2}$ is an eigenvalue of the magnetic Dirac operator $D^{t\eta}$. 
On the other hand, we compute for $n=2m+1$ and $b=-t$
$$\frac{n}{4(n-1)}\mathop{\rm inf}\limits_M\left(S-4t\lfloor\frac{n}{2}\rfloor^\frac{1}{2}|d\eta|\right)=\frac{2m+1}{8m}\left(2m(2m+4t+1)-8mt\right)=\frac{(2m+1)^2}{4}.
$$
Here, we use the fact that $|d\eta|^2=4m$.
Therefore, the equality in the magnetic Friedrich inequality is attained.\color{black} \\ 
\end{proof}

Simply-connected closed $\eta$-Einstein Sasakian spin manifolds can be constructed as circle bundles over simply-connected closed K\"ahler-Einstein manifolds with positive scalar curvature, see \cite[Example 6.1]{FriedrichKim:00}.
The scalar curvature of such a Sasaki manifold can be adjusted to the form above up to rescaling the metric on the K\"ahler-Einstein base.\\

In the following, we generalize B\"ar's \cite{Baer:92} and Hijazi's \cite{Hijazi:86} lower bounds for the smallest Dirac-eigenvalue to the magnetic Dirac operator, and, simultaneously, extend Reu\ss{}'s inequality \cite[Theorem 1]{Reuss:25} to nonzero eigenvalues of the magnetic Dirac operator.

\begin{theorem}\label{t:Hijazimagneticineq}
Let $\eta\in\Omega^1(M,\mathbb{R})$ be any real one-form on a closed connected spin manifold $(M^n,g)$.
Let $\lambda^\eta$ be any eigenvalue  of the magnetic Dirac operator $D^\eta$.
\begin{enumerate}
\item\label{item:Baermagnetic} If $n=2$ and the Euler characteristic $\chi(M)$ of $M$ is nonnegative, then 
\begin{equation}\label{eq:Baermagneticineq}
|\lambda^\eta|\geq\sqrt{\frac{2\pi\chi(M)}{\mathrm{Area}(M)}}-\|\eta\|_{L^\infty}.
\end{equation}
Moreover, equality is attained if and only if $\eta=0$ and $M$ is either a round $2$-sphere or a flat $2$-torus with trivial spin structure.
\item\label{item:Hijazimagnetic} If $n\geq3$ and $Y(M,[g])\geq0$, then 
\begin{equation}\label{eq:Hijazimagneticineq}
|\lambda^\eta|\mathrm{Vol}(M^n,g)^{\frac{1}{n}}\geq\sqrt{\frac{n}{4(n-1)}Y(M,[g])}-\|\eta\|_{L^n},
\end{equation}
where $\displaystyle Y(M,[g]):=\inf_{f\in C^\infty(M,\mathbb{R})\setminus\{0\}}\frac{\int_MfLf\,d\mu_g}{(\int_Mf^{\frac{2n}{n-2}}\,d\mu_g)^{\frac{n-2}{n}}}$ is the Yamabe invariant of $(M,[g])$ and $L:=4\frac{n-1}{n-2}\Delta+S$ is the Yamabe operator of $(M,g)$.
Moreover, equality when $\lambda^\eta\neq0$ and $\eta\neq0$ can only occur when $M$ is, up to rescaling the metric $g$, a Sasaki manifold admitting a real Killing spinor and with corresponding Reeb vector field $\eta$. 
If equality holds when $\lambda^\eta=0$ and $\eta\neq 0$, the manifold $M$ is conformally equivalent to an Einstein Sasaki manifold of positive scalar curvature whose Reeb vector field is the conformal change of $\eta$, which means that $(M,\overline{g}=e^{2u}g,\frac{e^{-u} \eta}{|\eta|})$ has to be Einstein-Sasaki of positive scalar curvature for some choice of the function $u$. \color{black}
\end{enumerate}
\end{theorem}

\begin{proof}
Let $\varphi\in\Gamma(\Sigma M)$. 
Using the Schr\"odinger-Lichnerowicz formula for $D$ as well as the pointwise identity $|\nabla\varphi|^2=|P\varphi|^2+\frac{1}{n}|D\varphi|^2$ involving the Penrose operator $P$, we have 
\begin{eqnarray*}
\int_M|P\varphi|^2\,d\mu_g&=&\int_M\left(|\nabla\varphi|^2-\frac{1}{n}|D\varphi|^2\right)d\mu_g\\
&=&\int_M\left(|D\varphi|^2-\frac{S}{4}|\varphi|^2-\frac{1}{n}|D\varphi|^2\right)d\mu_g\\
&=&\int_M\left(\frac{n-1}{n}|D\varphi|^2-\frac{S}{4}|\varphi|^2\right)d\mu_g,
\end{eqnarray*}
where $S$ is the scalar curvature of $(M^n,g)$.
It can be deduced that 
\[\int_M\left(|D\varphi|^2-\frac{n}{4(n-1)}S|\varphi|^2\right)d\mu_g=\frac{n}{n-1}\int_M|P\varphi|^2\,d\mu_g\geq0\]
for any spinor field $\varphi$ defined on $(M,g)$. Let $\psi\in\Gamma(\Sigma M)$ be any $\lambda^\eta$-eigenspinor for $D^{\eta}$ i.e., $D\psi=\lambda^\eta\psi-i\eta\cdot\psi$. If, instead of $g$, we consider any metric $\overline{g}:=e^{2u}g$ in its conformal class, where $u\in C^\infty(M,\mathbb{R})$ is arbitrary, then setting $\overline{\varphi}:=e^{-\frac{n-1}{2}u}\overline{\psi}$ and observing that $\overline{D}\overline{\varphi}=e^{-u}\left(\lambda^\eta\overline{\varphi}-i\overline{\eta}~\overline{\cdot}~\overline{\varphi}\right)$, we obtain 
\[\int_M\left(e^{-2u}|\lambda^\eta\overline{\varphi}-i\overline{\eta}~\overline{\cdot}~\overline{\varphi}|^2-\frac{n}{4(n-1)}\overline{S}|\overline{\varphi}|^2\right)d\mu_{\overline{g}}\geq0\]
that is,
\begin{equation}\label{eq:integineqHijazi}
\int_Me^{-u}|\lambda^\eta\psi-i\eta\cdot\psi|^2d\mu_g\geq\frac{n}{4(n-1)}\int_M\overline{S}e^{u}|\psi|^2d\mu_g.
\end{equation}
Here, we use the fact that $d\mu_{\overline{g}}=e^{nu}d\mu_g$. The key point is now to choose a suitable function $u$.\\

If $n=2$, then choosing $u\in C^\infty(M,\mathbb{R})$ such that $\displaystyle\Delta u=-\frac{S}{2}+\frac{\int_M S\,d\mu_g}{2\mathrm{Area}(M)}$, which is equivalent to $\displaystyle\Delta u=-\frac{S}{2}+\frac{2\pi\chi(M)}{\mathrm{Area}(M)}$ by the Gau\ss{}-Bonnet theorem, and using the identity $\overline{S}e^{2u}=S+2\Delta u$, the inequality \eqref{eq:integineqHijazi} becomes 
\[\int_Me^{-u}|\lambda^\eta\psi-i\eta\cdot\psi|^2d\mu_g\geq\frac{2\pi\chi(M)}{\mathrm{Area}(M)}\int_Me^{-u}|\psi|^2d\mu_g,\]
which yields 
\[(|\lambda^\eta|+\|\eta\|_{L^\infty})^2\geq\frac{2\pi\chi(M)}{\mathrm{Area}(M)}\]
and, provided $\chi(M)\geq0$, also \eqref{eq:Baermagneticineq}. 
Now, let us assume that equality holds in \eqref{eq:Baermagneticineq}.
Then $\overline{P}\overline{\varphi}=0$ i.e., $\overline{\varphi}$ must be a twistor-spinor on $(M,\overline{g}=e^{2u}g)$, the pointwise norm of $\eta$ must be constant and there must exist some nonnegative real number $\mu$ such that $i\eta\cdot\psi=-\mu\lambda^\eta\psi$.
If $\eta\neq0$, then necessarily $\mu>0$ and $\lambda^\eta\neq0$.
But the existence of a nowhere vanishing one-form $\eta$ on $M$ forces $\chi(M)=0$ and therefore $\lambda^\eta=0=|\eta|$ as well by the equality case, which is a contradiction.
Therefore $\eta=0$ on $M$ and the equality case is that of B\"ar's inequality, see \cite{Baer:92}.
This shows part (\ref{item:Baermagnetic}) of the Theorem.\\

When $n\geq3$, we first assume that $\psi_x\neq0$ for every $x\in M$.
Let $u:=\frac{2}{n-1}\ln(|\psi|)$ i.e., $\overline{g}=|\psi|^{\frac{4}{n-1}}g=f^{\frac{4}{n-2}}g$, where $f:=|\psi|^{\frac{n-2}{n-1}}$.
Observe then that $e^{-u}|\psi|^2=f^2$.
By $\overline{S}f^{\frac{4}{n-2}}=f^{-1}Lf$, the inequality \eqref{eq:integineqHijazi} becomes 
\[\int_M\big|\lambda^\eta\frac{\psi}{|\psi|}-i\eta\cdot\frac{\psi}{|\psi|}\big|^2f^2\,d\mu_g\geq\frac{n}{4(n-1)}\int_MfLf\,d\mu_g.\]
By H\"older inequality, we have 
\[\int_M\big|\lambda^\eta\frac{\psi}{|\psi|}-i\eta\cdot\frac{\psi}{|\psi|}\big|^2f^2\,d\mu_g\leq\big\|\lambda^\eta\frac{\psi}{|\psi|}-i\eta\cdot\frac{\psi}{|\psi|}\big\|_{L^n}^2\cdot\|f\|_{L^{\frac{2n}{n-2}}}^2,\]
from which 
\[\big\|\lambda^\eta\frac{\psi}{|\psi|}-i\eta\cdot\frac{\psi}{|\psi|}\big\|_{L^n}^2\geq\frac{n}{4(n-1)}\frac{\int_MfLf\,d\mu_g}{\|f\|_{L^{\frac{2n}{n-2}}}^2}\geq\frac{n}{4(n-1)}Y(M,[g])\]
and \eqref{eq:Hijazimagneticineq} follows in case $Y(M,[g])\geq0$.
As in \cite{FrankLoss:2022,Reuss:25}, let us now assume $\psi$ to possibly have zeroes on $M$ and set, for any $\varepsilon>0$, $|\psi|_\varepsilon:=\sqrt{|\psi|^2+\varepsilon^2}>0$.
Let $\overline{g}_\varepsilon:=e^{2u_\varepsilon}g=|\psi|_\varepsilon^{\frac{4}{n-1}}g=f_\varepsilon^{\frac{4}{n-2}}g$ as above, where $f_\varepsilon:=|\psi|_\varepsilon^{\frac{n-2}{n-1}}>0$.
We handle the left- and the right-hand-side of \eqref{eq:integineqHijazi} separately.
For the l.h.s., we can still write 
\[\int_Me^{-u_\varepsilon}\big|\lambda^\eta\psi-i\eta\cdot\psi\big|^2d\mu_g=\int_M\big|\lambda^\eta\frac{\psi}{|\psi|_\varepsilon}-i\eta\cdot\frac{\psi}{|\psi|_\varepsilon}\big|^2f_\varepsilon^2\,d\mu_g,\]
which is bounded above by $(|\lambda^\eta|\mathrm{Vol}(M,g)^{\frac{1}{n}}+\|\eta\|_{L^n})^2\cdot\|f_\varepsilon\|_{L^{\frac{2n}{n-2}}}^2$ taking $|\psi|<|\psi|_\varepsilon$ into account.
As for the r.h.s. of \eqref{eq:integineqHijazi}, we have 
\begin{eqnarray*}
\int_M\overline{S}e^{u_\varepsilon}|\psi|^2d\mu_g&=&\int_M\overline{S}e^{2u_\varepsilon}e^{-u_\varepsilon}|\psi|_\varepsilon^2\cdot\frac{|\psi|^2}{|\psi|_\varepsilon^2}d\mu_g\\
&=&\int_M\overline{S}f_\varepsilon^{\frac{4}{n-2}}\cdot f_\varepsilon^2\cdot\frac{|\psi|^2}{|\psi|_\varepsilon^2}d\mu_g\\
&=&\int_Mf_\varepsilon Lf_\varepsilon\cdot \frac{|\psi|^2}{|\psi|_\varepsilon^2}d\mu_g.
\end{eqnarray*}
Now $\displaystyle\frac{|\psi|^2}{|\psi|_\varepsilon^2}=1-\frac{\varepsilon^2}{|\psi|_\varepsilon^2}=1-\varepsilon^2f_\varepsilon^{-\frac{2(n-1)}{n-2}}$, so that 
\begin{eqnarray}\label{eq:estimLfepsilonpsi}
\nonumber\int_Mf_\varepsilon Lf_\varepsilon\cdot \frac{|\psi|^2}{|\psi|_\varepsilon^2}d\mu_g&=&\frac{4(n-1)}{n-2}\int_M\langle df_\varepsilon,d(f_\varepsilon\cdot\frac{|\psi|^2}{|\psi|_\varepsilon^2})\rangle\,d\mu_g+\int_MSf_\varepsilon^2\frac{|\psi|^2}{|\psi|_\varepsilon^2}\,d\mu_g\\
\nonumber&=&\frac{4(n-1)}{n-2}\int_M\langle df_\varepsilon,d(f_\varepsilon-\varepsilon^2f_\varepsilon^{-\frac{n}{n-2}})\rangle\,d\mu_g\\
\nonumber&&+\int_MSf_\varepsilon^2(1-\varepsilon^2f_\varepsilon^{-\frac{2(n-1)}{n-2}})\,d\mu_g\\
\nonumber&=&\frac{4(n-1)}{n-2}\int_M(1+\frac{n\varepsilon^2}{n-2}f_\varepsilon^{-\frac{2(n-1)}{n-2}})|df_\varepsilon|^2\,d\mu_g\\
\nonumber&&+\int_MSf_\varepsilon^2(1-\varepsilon^2f_\varepsilon^{-\frac{2(n-1)}{n-2}})\,d\mu_g\\
\nonumber&=&\int_Mf_\varepsilon Lf_\varepsilon\,d\mu_g
+\frac{4n(n-1)\varepsilon^2}{(n-2)^2}\int_Mf_\varepsilon^{-\frac{2(n-1)}{n-2}}|df_\varepsilon|^2\,d\mu_g\\
\nonumber&&-\varepsilon^2\int_M Sf_\varepsilon^{-\frac{2}{n-2}}\,d\mu_g\\
&\geq&\|f_\varepsilon\|_{L^{\frac{2n}{n-2}}}^2\cdot\left(Y(M,[g])-\varepsilon^2\frac{\int_M Sf_\varepsilon^{-\frac{2}{n-2}}\,d\mu_g}{\|f_\varepsilon\|_{L^{\frac{2n}{n-2}}}^2}\right).
\end{eqnarray}
Therefore, we get
\begin{equation}\label{eq:Hijazialmostdone}(|\lambda^\eta|\mathrm{Vol}(M,g)^{\frac{1}{n}}+\|\eta\|_{L^n})^2\geq\frac{n}{4(n-1)}\Big(Y(M,[g])-\varepsilon^2\frac{\int_M Sf_\varepsilon^{-\frac{2}{n-2}}\,d\mu_g}{\|f_\varepsilon\|_{L^{\frac{2n}{n-2}}}^2}\Big).\end{equation}
Next, we inspect the only term involving $\varepsilon$ as $\varepsilon\to0$.
Obviously, we have $\|f_\varepsilon\|_{L^{\frac{2n}{n-2}}}\underset{\varepsilon\to0}{\longrightarrow}\|f\|_{L^{\frac{2n}{n-2}}}>0$ with $f:=|\psi|^{\frac{n-2}{n-1}}\geq0$ and where $\psi$ does not vanish identically.
As for $\varepsilon^2f_\varepsilon^{-\frac{2}{n-2}}$, it is sufficient to observe that $f_\varepsilon\geq\varepsilon^{\frac{n-2}{n-1}}$ to obtain 
\[\varepsilon^2f_\varepsilon^{-\frac{2}{n-2}}\leq\varepsilon^{2-\frac{2}{n-1}}=\varepsilon^{\frac{2(n-2)}{n-1}}\]
uniformly on $M$, so that $\displaystyle\varepsilon^2\frac{\int_M Sf_\varepsilon^{-\frac{2}{n-2}}\,d\mu_g}{\|f_\varepsilon\|_{L^{\frac{2n}{n-2}}}^2}\underset{\varepsilon\to0}{\longrightarrow}0$ and the inequality \eqref{eq:Hijazimagneticineq} follows from \eqref{eq:Hijazialmostdone} in the case of $Y(M,[g])\geq0$. Assume now \eqref{eq:Hijazimagneticineq} to be an equality.
We mimic the proof of \cite[Theorem 22]{Reuss:25}.
Fix $\varepsilon>0$.
Then, for the eigenspinor $\psi$ and the function $f_\varepsilon:=|\psi|^{\frac{n-2}{n-1}}\geq0$ defined above, we have by \eqref{eq:estimLfepsilonpsi} that
\begin{equation}\label{eq:equalcaseHijazi}
(|\lambda^\eta|\mathrm{Vol}(M,g)^{\frac{1}{n}}+\|\eta\|_{L^n})^2\cdot\|f_\varepsilon\|_{L^{\frac{2n}{n-2}}}^2
\geq\frac{n}{4(n-1)}\Big(\int_Mf_\varepsilon Lf_\varepsilon\;d\mu_g+A(\varepsilon)\Big),
\end{equation}
where $A(\varepsilon):=\frac{4n(n-1)\varepsilon^2}{(n-2)^2}\int_Mf_\varepsilon^{-\frac{2(n-1)}{n-2}}|df_\varepsilon|^2\,d\mu_g-\varepsilon^2\int_M Sf_\varepsilon^{-\frac{2}{n-2}}\,d\mu_g$.
Now, replacing the expression $(|\lambda^\eta|\mathrm{Vol}(M,g)^{\frac{1}{n}}+\|\eta\|_{L^n})^2$ by $\frac{n}{4(n-1)}Y(M,[g])$, we observe that 
\[Y(M,[g])\geq\frac{\int_Mf_\varepsilon Lf_\varepsilon\;d\mu_g}{\|f_\varepsilon\|_{L^{\frac{2n}{n-2}}}^2
}+\frac{A(\varepsilon)}{\|f_\varepsilon\|_{L^{\frac{2n}{n-2}}}^2}\geq Y(M,[g])+\frac{A(\varepsilon)}{\|f_\varepsilon\|_{L^{\frac{2n}{n-2}}}^2},\]
from which $A(\varepsilon)\leq0$ can be deduced.
But, because of $\varepsilon^2\int_M Sf_\varepsilon^{-\frac{2}{n-2}}\,d\mu_g\underset{\varepsilon\to0}{\longrightarrow}0$ as we have seen above, 
\[\varepsilon^2\int_Mf_\varepsilon^{-\frac{2(n-1)}{n-2}}|df_\varepsilon|^2\,d\mu_g\underset{\varepsilon\to0}{\longrightarrow}0\]
must hold.
Note that we have used the inequality $\int_M Sf_\varepsilon^{-\frac{2}{n-2}}\,d\mu_g\geq0$ which is a consequence of  \eqref{eq:Hijazialmostdone}.
Now, since $A(\varepsilon)\underset{\varepsilon\to0}{\longrightarrow}0$ and $f_\varepsilon\underset{\varepsilon\to0}{\longrightarrow}f=|\psi|^{\frac{n-2}{n-1}}\neq0$ in any $L^p$-norm, $p>1$, the $L^2$-norm of $df_\varepsilon$ must remain bounded as $\varepsilon\to0$.
Namely, \eqref{eq:equalcaseHijazi} implies that 
\begin{eqnarray*}
\frac{n}{n-2}\|df_\varepsilon\|_{L^2}^2&\leq&(|\lambda^\eta|\mathrm{Vol}(M,g)^{\frac{1}{n}}+\|\eta\|_{L^n})^2\cdot\|f_\varepsilon\|_{L^{\frac{2n}{n-2}}}^2\\
&&-\frac{n}{4(n-1)}\Big(\int_MSf_\varepsilon^2\,d\mu_g+A(\varepsilon)\Big),
\end{eqnarray*}
whose r.h.s. remains bounded as $\varepsilon\to0$.
The key argument is now \cite[Lemma 9]{FrankLoss:2022}, which states that, under our assumptions, the function $f$ must be weakly differentiable and that $\|df\|_{L^2}^2\leq\underset{\varepsilon\to0}{\liminf}\|df_\varepsilon\|_{L^2}^2$.
Therefore, as $\varepsilon\to0$, we obtain from \eqref{eq:equalcaseHijazi} that 
\[Y(M,[g])\geq\frac{\int_M\frac{4(n-1)}{n-2}|df|^2+Sf^2\,d\mu_g}{\|f\|_{L^{\frac{2n}{n-2}}}^2}\geq Y(M,[g]),\]
such that $Y(M,[g])=\frac{\int_M\frac{4(n-1)}{n-2}|df|^2+Sf^2\,d\mu_g}{\|f\|_{L^{\frac{2n}{n-2}}}^2}$ actually holds, which shows $f$ to minimize the Yamabe functional and therefore to be \emph{positive}.
This implies that $\psi$ has no zero on $M$.
Therefore we can use the first part of the proof of inequality \eqref{eq:Hijazimagneticineq} where the eigenspinor $\psi$ is assumed to have no zero on $M$.\\

We first handle the case where $\lambda^\eta\neq0$.
Since \eqref{eq:Hijazimagneticineq} is an equality, $\overline{P}\overline{\varphi}=0$ must first hold for the conformal metric $\overline{g}=|\psi|^{\frac{4}{n-1}}g=f^{\frac{4}{n-2}}g$ on $M$ i.e., $\overline{\varphi}=\frac{\overline{\psi}}{|\psi|}$ is a twistor-spinor on $(M,\overline{g})$.
Observe that, since $\overline{\varphi}$ is a twistor-spinor of \emph{constant length} $1$ on $(M,\overline{g})$, it must be either a parallel spinor or the sum of two real (nonparallel) Killing spinors for opposite Killing constants.
Actually $\psi$ must be a parallel or real Killing spinor.
Namely, equality in the Minkowski inequality $\|\lambda^\eta\frac{\psi}{|\psi|}-i\eta\cdot\frac{\psi}{|\psi|}\|_{L^n}\leq|\lambda^\eta|\mathrm{Vol}(M,g)^{\frac{1}{n}}+\|\eta\|_{L^n}$ implies the existence of a nonnegative real number $\mu$ on $M$ such that $i\eta\cdot\psi=-\mu\lambda^\eta\psi$ on $M$.
If $\eta$ vanishes at some point, then $\mu$ or $\lambda^\eta$ vanishes as well and thus $\eta=0$ on $M$, such that the equality case is just that of Hijazi's inequality for the standard Dirac operator.
Otherwise, necessarily $\mu>0$, the function $|\eta|=\mu|\lambda^\eta|>0$ is constant on $M$ and 
\begin{equation}\label{eq:Cliffactioneta}
i\frac{\eta}{|\eta|}\cdot\psi=\epsilon\psi
\end{equation}
must hold on $M$, where $\epsilon=-\mathrm{sgn}(\lambda^\eta)\in\{\pm1\}$.
In turn, this implies that $i\overline{\eta}~\overline{\cdot}~\overline{\varphi}=\epsilon|\eta|\overline{\varphi}=\epsilon\mu|\lambda^\eta|\overline{\varphi}=-\mu\lambda^\eta\overline{\varphi}$ holds as well.
But then $\overline{D}\overline{\varphi}=e^{-u}\left(\lambda^\eta\overline{\varphi}-i\overline{\eta}~\overline{\cdot}~\overline{\varphi}\right)=\lambda^\eta(1+\mu)e^{-u}\overline{\varphi}$, which implies with $\overline{P}\overline{\varphi}=0$ that 
\[\overline{\nabla}_X\overline{\varphi}=-\frac{1}{n}X\overline{\cdot}~\overline{D}~\overline{\varphi}=-\frac{\lambda^\eta(1+\mu)e^{-u}}{n}X~\overline{\cdot}
~\overline{\varphi}\]
for all tangent vectors $X$ to $M$.
By \cite[Cor. 3.6]{Hijazi:86}, such an identity implies the function $u$ to be constant and therefore $\overline{\varphi}$ to be a real Killing spinor on $(M,\overline{g})$.
Because $u$ (and hence $|\psi|$) is constant on $M$, the spinor $\psi$ must be a real Killing spinor on $(M,g)$.
Up to rescaling $\psi$, we may assume that $|\psi|=1$ on $M$.
It remains to note that, up to rescaling $g$ so as to make $|\eta|=1$, the identity \eqref{eq:Cliffactioneta} provides a Sasaki structure on $(M,g)$, in particular $M$ must be odd-dimensional.
This must be long known, but we write down an elementary proof for the sake of being self-contained.
Up to replacing $\eta$ by $\epsilon\eta$, it can be assumed that $\epsilon=1$ in \eqref{eq:Cliffactioneta}.
Note that $\eta$ is a unit Killing vector field then because it coincides (up to sign) with the Killing vector field which is naturally associated with $\psi$.
Up to rescaling the metric $g$ -- and $\eta$ accordingly -- on $M$, it can be assumed that $\psi$ is an $\frac{\epsilon}{2}$-Killing spinor on $(M,g)$ for some $\epsilon\in\{\pm1\}$.
Then, for every tangent vector $X$ on $M$, we have $i\nabla_X\eta\cdot\psi+i\eta\cdot\nabla_X\psi=\epsilon\nabla_X\psi$ that is,
\begin{equation}\label{eq:nablaXetacdotpsi}\nabla_X\eta\cdot\psi=\epsilon(g(X,\eta)\psi-iX\cdot\psi)\end{equation}
on $M$.
Replacing $X$ by $\nabla_X\eta$, it follows that, for all $X\in\eta^\perp$,
\[\nabla_{\nabla_X\eta}\eta\cdot\psi=-X\cdot\psi,\]
meaning that $(\nabla\eta)^2=-\mathrm{Id}_{\eta^\perp}$ holds on $\eta^\perp\subset TM$.
It remains to notice that $\nabla_X(\nabla\eta)(Y)=g(Y,\eta)X-g(X,Y)\eta$ holds for all $X,Y\in TM$, which follows from differentiating \eqref{eq:nablaXetacdotpsi} again.
Namely, for all tangent vectors $X,Y$ on $M$,
\begin{eqnarray*}
\nabla_X\nabla_Y\eta\cdot\psi&=&\nabla_X\left(\nabla_Y\eta\cdot\psi\right)-\nabla_Y\eta\cdot\nabla_X\psi\\
&=&\epsilon\left(g(\nabla_XY,\eta)+g(Y,\nabla_X\eta)\right)\psi-i\epsilon\nabla_XY\cdot\psi\\
&&+\epsilon\left(g(Y,\eta)-iY\cdot\right)\nabla_X\psi-\frac{\epsilon}{2}\nabla_Y\eta\cdot X\cdot\psi\\
&=&\epsilon g(\nabla_XY,\eta)\psi-i\epsilon\nabla_XY\cdot\psi+\epsilon\Big(\underbrace{g(Y,\nabla_X\eta)+g(X,\nabla_Y\eta)}_{0}\Big)\psi\\
&&+\frac{1}{2}g(Y,\eta)X\cdot\psi-\frac{i}{2}Y\cdot X\cdot\psi+\frac{\epsilon}{2}X\cdot\nabla_Y\eta\cdot\psi\\
&=&\epsilon g(\nabla_XY,\eta)\psi-i\epsilon\nabla_XY\cdot\psi+g(Y,\eta)X\cdot\psi-\frac{i}{2}(Y\cdot X+X\cdot Y)\cdot\psi\\
&=&\epsilon g(\nabla_XY,\eta)\psi-i\epsilon\nabla_XY\cdot\psi+g(Y,\eta)X\cdot\psi+ig(X,Y)\psi,
\end{eqnarray*}
so that, recalling that $i\eta\cdot\psi=\psi$, we obtain
\begin{eqnarray*}
\nabla_X(\nabla\eta)(Y)\cdot\psi&=&\nabla_X\nabla_Y\eta\cdot\psi-\nabla_{\nabla_XY}\eta\cdot\psi\\
&=&\epsilon g(\nabla_XY,\eta)\psi-i\epsilon\nabla_XY\cdot\psi+g(Y,\eta)X\cdot\psi+ig(X,Y)\psi\\
&&-\epsilon(g(\nabla_XY,\eta)\psi-i\nabla_XY\cdot\psi)\\
&=&g(Y,\eta)X\cdot\psi+ig(X,Y)\psi\\
&=&g(Y,\eta)X\cdot\psi-g(X,Y)\eta\cdot\psi.
\end{eqnarray*}
Therefore, $(M,g,\eta)$ must be a Sasaki manifold.\\

Assume now equality is realized for $\lambda^\eta=0$, we then have as before that $\overline\varphi$ is a twistor spinor of length $1$. The manifold $(M,\overline{g})$ has to be Einstein of scalar curvature ${\rm vol}(M^n,\overline{g})^{-\frac{2}{n}}Y(M,[g])$. Namely, the scalar curvature of $\overline{g}$ has to be equal to 
\begin{equation}\label{eq:scalconf}
\overline{S}=\frac{4(n-1)}{n}\frac{|\overline{D}\overline{\varphi}|^2}{|\overline\varphi|^2}=\frac{4(n-1)}{n}e^{-2u}|\eta|^2={\rm const},
\end{equation}
since we have $\overline{D}\overline{\varphi}=-ie^{-u}\overline{\eta}~\overline{\cdot}~\overline{\varphi}$. From the fact that $\|\eta\|_{L^n}=\sqrt{\frac{n}{4(n-1)}Y(M,[g])}$, we deduce that this constant must be ${\rm vol}(M^n,\overline{g})^{-\frac{2}{n}}Y(M,[g])$. If $\overline S$ is zero, then $\eta$ must be zero which is a contradiction. 
Hence $\overline S$ is positive, meaning that $\eta$ cannot have any zero. 
In this case, $\overline \varphi$ is a sum of two Killing spinors of different constants. 
Indeed, $\overline\varphi$ can be written as $\overline{\varphi}=\varphi_++\varphi_-$ where $$\varphi_\pm=\frac{1}{2}(\overline\varphi\pm\frac{1}{c}\overline D\overline \varphi)=\frac{1}{2}(\overline\varphi\mp\frac{i}{|\eta|}\overline\eta~\overline{\cdot}~\overline \varphi).$$ 
Here, we have $c:=\sqrt{\frac{n\overline S}{4(n-1)}}=e^{-u}|\eta|$ and $\overline{D}\overline{\varphi}=-ie^{-u}\overline{\eta}~\overline{\cdot}~\overline{\varphi}$.
The spinors $\varphi_\pm$ are Killing spinors associated with the Killing numbers $\mp \frac{e^{-u}|\eta|}{n}$. 
Moreover, an easy computation shows that 
$$\overline\eta~\overline\cdot~\varphi_\pm=\pm i|\eta|\varphi_\pm.$$
Hence, we are in the same case as before. 
Therefore, $(M,\overline{g},\frac{\overline \eta}{|\eta|})$ has to be Einstein-Sasaki of positive scalar curvature. 
\end{proof}

Note that two further inequalities can be deduced from \eqref{eq:integineqHijazi} in case $n\geq3$: on the one hand, 
\[\displaystyle(|\lambda^\eta|+\|\eta\|_{L^\infty})^2\geq\frac{n}{4(n-1)}\sup_{u\in C^\infty(M,\mathbb{R})}\inf_M(\overline{S}e^{2u}),\]
on the other hand 
\[\displaystyle(|\lambda^\eta|+\|\eta\|_{L^\infty})^2\geq\frac{n}{4(n-1)}\mu_1,\]
where $\mu_1$ is the smallest eigenvalue of $L$.
Those inequalities, along with \eqref{eq:Hijazimagneticineq}, extend Hijazi's inequalities \cite{Hijazi:86} to the magnetic Dirac operator.
Remark in particular that \eqref{eq:Hijazimagneticineq} is a conformal lower bound for the normalized quantity $|\lambda^\eta|\mathrm{Vol}(M^n,g)^{\frac{1}{n}}$.\\

\begin{remark}\label{rem:proofReuss}\rm We give a simple proof of Reu\ss{}'s result \cite[Theorem 1]{Reuss:25}.
Let $\psi\neq0$ be a spinor satisfying $D^\eta\psi=0$. 
Assume first that $\eta$ has no zeroes and let us make the conformal change of the metric $\overline{g}:=|\eta|^2 g=e^{2u}g$. 
As above, we let $\overline\varphi:=e^{-\frac{n-1}{2}u}\overline\psi$ and obtain
$$\overline{D}\overline\varphi=-ie^{-u}\overline\eta~\overline\cdot~\overline\varphi
=-i|\eta|^{-1}\overline\eta~\overline\cdot~\overline\varphi.$$

Therefore, by the min-max principle, we find that 
$$(\lambda^{\eta,{\overline g}})^2\leq \frac{\int_M|\overline D\overline \varphi|_{\overline g}^2 d\mu_{\overline g}}{\int_M |\overline \varphi|^2 d\mu_{\overline g}}=\frac{\int_M|\eta|^{-2}|\overline\eta|_{\overline g}^2|\overline \varphi|^2 d\mu_{\overline g}}{\int_M |\overline \varphi|^2 d\mu_{\overline g}}=1.$$
Thus, by the standard Hijazi inequality, we deduce that 
$$\sqrt{\frac{n}{4(n-1)}Y(M,[g])}\leq |\lambda^{\eta,{\overline{g}}}|\mathrm{Vol}(M^n,\overline{g})^{\frac{1}{n}}\leq\mathrm{Vol}(M^n,\overline{g})^{\frac{1}{n}}=\|\eta\|_{L^n},$$ 
which is the inequality of \cite[Theorem 1]{Reuss:25}.
Note that this still works in dimension $n=2$ up to replacing $\frac{n}{4(n-1)}Y(M,[g])$ by $2\pi\chi(M)$. 
If now $\eta$ does have zeroes on $M$, then as before we let $|\eta|_\epsilon:=\sqrt{(|\eta|^2+\epsilon^2)}$ for any $\epsilon>0$ and set $\overline{g}_\epsilon:=|\eta|_\epsilon^2g$ on $M$.
If $D^\eta\psi=0$ on $M$, then we let $\overline{\varphi}:=|\eta|_\epsilon^{-\frac{n-1}{2}}\overline{\psi}$.
As above, $\overline{D}\overline{\varphi}=-|\eta|_\epsilon^{-1}i\overline{\eta}~\overline{\cdot}~\overline{\varphi}$ holds and, by the min-max principle, we obtain 
\begin{eqnarray*}
\lambda^\eta_1(\overline{D})^2&\leq&\frac{\int_M|\overline{D}\overline{\varphi}|^2d\mu_{\overline{g}_\epsilon}}{\int_M|\overline{\varphi}|^2d\mu_{\overline{g}_\epsilon}}\\
&=&\frac{\int_M\frac{|\eta|^2}{|\eta|_\epsilon^2}|\overline{\varphi}|^2d\mu_{\overline{g}_\epsilon}}{\int_M|\overline{\varphi}|^2d\mu_{\overline{g}_\epsilon}}\\
&=&1-\epsilon^2\frac{\int_M\frac{|\overline{\varphi}|^2}{|\eta|_\epsilon^2}d\mu_{\overline{g}_\epsilon}}{\int_M|\overline{\varphi}|^2d\mu_{\overline{g}_\epsilon}}\\
&\leq&1,
\end{eqnarray*}
so that, by B\"ar's resp. Hijazi's inequality in dimension $n=2$ resp. $n\geq3$, we obtain 
\begin{equation}\label{eq:minmaxgepsilon|eta|}
\||\eta|_\epsilon\|_{L^n}^2=\mathrm{Vol}(M,\overline{g}_\epsilon)^{\frac{2}{n}}\geq\lambda^\eta_1(\overline{D})^2\mathrm{Vol}(M,\overline{g}_\epsilon)^{\frac{2}{n}}\geq\left\{\begin{array}{ll}2\pi\chi(M)&\textrm{for }n=2\\\frac{n}{4(n-1)}Y(M,[g])&\textrm{for }n\geq3.
\end{array}\right.
\end{equation}
It remains to notice that, since $\left\||\eta|_\epsilon-|\eta|\right\|_{L^\infty}\underset{\epsilon\to0}{\longrightarrow}0$, we have $\||\eta|_\epsilon\|_{L^n}^2\underset{\epsilon\to0}{\longrightarrow}\|\eta\|_{L^n}^2$ and the same inequality follows.
Note that this enhances inequality \eqref{eq:Baermagneticineq} in dimension $n=2$ and for $\lambda^\eta=0$ since we now have $\|\eta\|_{L^2}$ instead of $\|\eta\|_{L^\infty}$.
\end{remark}

It is essential here to point out that inequality 
\[\left(|\lambda^\eta|+\|\eta\|_{L^\infty}\right)^2\geq\frac{2\pi\chi(M)}{\mathrm{Area}(M)},\]
which is a reformulation of \eqref{eq:Baermagneticineq} without the square root, is actually trivial when the genus of $M$ is positive (including for the $2$-torus), because then $\chi(M)\leq0$.
When $M=\mathbb{S}^2$, it turns out that $0$ cannot be a $D^\eta$-eigenvalue. 
To see that, we consider the general case of a K\"ahler spin manifold $(M^{2m},g,J)$. 
Recall that the action of the K\"ahler form $\Omega(\cdot,\cdot)=g(J\cdot,\cdot)$ splits the spinor bundle $\Sigma M$ into $\displaystyle\Sigma M=\bigoplus_{r=0}^m\Sigma_r M$, where $\Sigma_rM$ is the eigenspace associated with the eigenvalue $i(2r-m)$ of $\Omega$ for each $r\in\{0,\ldots,m\}$. 
Moreover, for any $X\in TM$, $p_-(X)\cdot\Sigma_r M\in\Sigma_{r-1} M$  and $p_+(X)\cdot\Sigma_r M\in\Sigma_{r+1} M$, where $p_\pm(X):=\frac{1}{2}(X\mp iJX)$ is the standard projector onto pointwise $\pm i$-eigenspaces of $J$. 
Also, we have that $\displaystyle\Sigma_+M=\bigoplus_{r\,\, \rm even}\Sigma_r M$ and $\displaystyle\Sigma_-M=\bigoplus_{r\,\, \rm odd}\Sigma_r M$.

\begin{proposition}\label{p:DetaKaehler}
Let $(M^{2m},g,J)$ be a K\"ahler spin manifold. 
Then for any $\eta\in \Omega^1(M)$ and smooth function $f$, we have
\begin{equation}\label{eq:invariancekahler}
D^{\eta+Jdf}=e^fD^\eta e^{-f}
\end{equation}
on $\Sigma_0 M$. 
The same identity holds on $\Sigma_m M$ by replacing $f$ by $-f$ on the r.h.s.
In particular, on $\mathbb{S}^2$, we have $\ker(D^\eta)=\ker(D)=\{0\}$. 
\end{proposition}
\begin{proof} For any section $\psi_0\in \Sigma_0 M$ we have
\begin{eqnarray*}
D^{\eta+Jdf}\psi_0&=&D^\eta\psi_0+iJdf\cdot\psi_0\;\textrm{where}\; D\psi_0\in\Gamma(\Sigma_1M)\\
&=&D^\eta\psi_0+ip_+(Jdf)\cdot\psi_0\;\textrm{since}\;p_-(Jdf)\cdot\psi_0=0\\
&=&D^\eta\psi_0-p_+(df)\cdot\psi_0\;\textrm{because}\;p_+(Jdf)=ip_+(df)\\
&=&D^\eta\psi_0-df\cdot\psi_0\;\textrm{since}\;p_-(df)\cdot\psi_0=0\\
&=&e^fD^\eta(e^{-f}\psi_0).
\end{eqnarray*}
The computation of $D^{\eta+Jdf}\psi_m$ is similar. 
To prove the second part, we use the Hodge decomposition theorem.
Namely, there exists $h\in C^\infty(M,\mathbb{R})$ as well as $\omega\in\Omega^2(M)$ such that $\eta=dh+\delta\omega$ (there is no nonzero harmonic $1$-form on $\mathbb{S}^2$). 
By Proposition \ref{prop:gaugeonva}, since $D^{dh+\delta\omega}$ and $D^{\delta\omega}$ are unitarily equivalent, it is sufficient to assume that $\eta=\delta\omega$ is co-exact, hence $\eta=\delta(f\mathrm{vol}_g)=-\nabla f\lrcorner\,\mathrm{vol}_g=-Jdf$ for some real-valued function $f$, where $J$ is the natural K\"ahler structure on $M$ provided by its orientation. Let $\psi\in {\rm ker}(D^\eta)$. 
By writing $\psi=\psi_0+\psi_1\in\Gamma(\Sigma_0 M\oplus\Sigma_1M)$, we get $\psi_0$ and $\psi_1\in {\rm ker}(D^{-Jdf})$. 
Hence, from \eqref{eq:invariancekahler}, we deduce $e^f\psi_0$ and $e^{-f}\psi_1$ are in ${\rm ker}(D)=\{0\}$. 
Thus, $\psi_0=\psi_1=0$ and therefore, $\psi=0$. 
\end{proof}
\color{black}

Apart from generalizing Friedrich's inequality we also present an estimate
for the eigenvalues of the magnetic Dirac operator that involves the nodal set of the corresponding eigenspinors which is
based on the ideas presented in \cite{Branding:2018a,Branding:2018b}. 
More precisely, we have the following

\begin{theorem}\label{t:estimnodalset}
Let \((M,g)\) be a closed spin surface
and \(\lambda^2_k\) be the \(k\)-th eigenvalue of the magnetic Dirac operator \((D^\eta)^2\).
Then the following eigenvalue estimate holds
\begin{equation}
\label{eigenvalue-estimate}
\lambda_k^2\geq\frac{2\pi\chi(M)}{\vol(M,g)}
-\frac{\int_M|d\eta|d\mu_g}{\vol(M,g)}
+\frac{4\pi N_k}{\vol(M,g)},
\end{equation}
where \(\chi(M)\) is the Euler characteristic of \(M\) 
and \(N_k\) denotes the sum of the order of the zeros of an eigenspinor \(\psi_k\)
belonging to the \(k\)-th eigenvalue of the magnetic Dirac operator, 
that is
\begin{equation}
\label{nodal-set}
N_k=\max\left(\sum_{p\in M,|\psi_k|(p)=0}n_p\right).
\end{equation}

In the case that \(\psi\in\Gamma(\Sigma^\pm M)\) is in the kernel of the magnetic Dirac operator we have
\begin{align*}
N_0(\psi)=-\frac{\chi(M)}{2}.   
\end{align*}
\end{theorem}

\begin{proof}
By the main result of \cite{Baerzerosets:99} we know that on a two-dimensional manifold the zero-set of any
eigenspinor of the magnetic Dirac operator is discrete. 
In the following we will make use of the energy-momentum tensor \(T^\eta(X,Y)\)
associated with the magnetic Dirac operator, which is defined as follows
\[
T^\eta(X,Y):=\langle X\cdot\nabla^\eta_Y\psi+Y\cdot\nabla^\eta_X\psi,\psi\rangle,
\]
where \(X,Y\in\Gamma(TM)\).

The following equation is a version of \cite[Lemma 4.2]{FriedrichKim:01} 
adapted to the case of the magnetic Dirac operator, i.e.
\begin{equation}
\label{inequality-modified-connection}
\frac{\langle\psi,(D^\eta)^2\psi\rangle}{|\psi|^2}\geq\frac{S}{4}+\frac{1}{2}\frac{\langle id\eta\cdot\psi,\psi\rangle}{|\psi|^2}+\frac{|T^\eta|^2}{4|\psi|^4}+\Delta\log|\psi|
-\frac{\langle D^\eta\psi,d(\log|\psi|^2)\cdot\psi\rangle}{|\psi|^2},
\end{equation}
which holds away from the zero-set of \(\psi\), see also \cite[Lemma 2.1]{Branding:2018b}. 
Note that a different sign convention for the Laplace operator was used in \cite{Branding:2018a,Branding:2018b}.
In order to establish \eqref{inequality-modified-connection}
we define a new magnetic connection on the spinor bundle by
\begin{align*}
    \tilde\nabla^\eta_X\psi:=\nabla^\eta_X\psi-2\alpha(X)\psi-\beta(X)\cdot\psi-X\cdot\alpha\cdot\psi
\end{align*}
with a one-form \(\alpha\) and a symmetric \((1,1)\)-tensor \(\beta\) which are given by
\begin{align*}
    \alpha:=\frac{d|\psi|^2}{2|\psi|^2},\qquad \beta:=-\frac{T^\eta(\cdot,\cdot)}{2|\psi|^2}.
\end{align*}
A direct calculation then shows that
\begin{align*}
0\leq|\tilde\nabla^\eta\psi|^2=&|\nabla^\eta\psi|^2+2|\alpha|^2|\psi|^2+|\beta|^2|\psi|^2    
-4\alpha(e_i)\langle\nabla^\eta_{e_i}\psi,\psi\rangle \\
&+2\langle\psi,\beta(e_i)\cdot\nabla^\eta_{e_i}\psi\rangle
+2\langle D^\eta\psi,\alpha\cdot\psi\rangle.
\end{align*}
Now, direct calculations using \eqref{eq:SchroedingerLichnerowiczmagnetic} show that
\begin{align*}
    |\nabla^\eta\psi|^2=&-\Delta\frac{1}{2}|\psi|^2+\langle\psi,(D^\eta)^2\psi\rangle-\frac{S}{4}|\psi|^2
    -\langle\psi,id\eta\cdot\psi\rangle
    ,\\
    \alpha(e_i)\langle\nabla^\eta_{e_i}\psi,\psi\rangle=&|\alpha|^2|\psi|^2=\frac{\big|d|\psi|^2\big|^2}{4|\psi|^2}, \\
    \langle\beta(e_i)\cdot\nabla^\eta_{e_i}\psi,\psi\rangle=&-\frac{|T^\eta|^2}{4|\psi|^2}.
\end{align*}

Moreover, recall that on a closed Riemannian surface we have that if the zero set of \(|\psi|\) is discrete and \(|\psi|\) does not vanish identically, then the following equality holds 
\begin{equation}
\label{laplace-log}
\mathrm{p.v.}\left(\Delta\log|\psi|\right)=2\pi\sum_{p\in M,|\psi|(p)=0}n_p,
\end{equation}
where $\mathrm{p.v.}$ stands for the Cauchy principal value of the distribution and \(n_p\) is the order of \(|\psi|\) at the point \(p\), see \cite{Branding:2018a} for a proof.
As a next step we apply \eqref{inequality-modified-connection} for \(\psi\) being an eigenspinor of \((D^\eta)^2\)
in which case we can also estimate the energy-momentum tensor by \(|T^\eta|^2\geq 2\lambda^2|\psi|^4\).
Thus, from \eqref{inequality-modified-connection} we obtain
\[
\lambda^2\geq\frac{K}{2}-|d\eta|+2\Delta\log|\psi|,
\]
where \(K\) denotes the Gaussian curvature of \(M\). 
The first claim now follows from integrating the above inequality and application of the Gauss-Bonnet theorem.

Regarding the second claim on the zero-set of spinors in the kernel of the magnetic Dirac operator we 
assume that \(\psi\in\Gamma(\Sigma^+M)\) is in the kernel of the magnetic Dirac operator \(D^\eta\) such that we get
\begin{align*}
-\Delta\log|\psi|=\frac{K}{2}
+\frac{1}{2}\frac{\langle id\eta\cdot\psi,\psi\rangle}{|\psi|^2}+\frac{|\nabla^\eta\psi|^2}{|\psi|^2}
-\frac{1}{2}\frac{\big|d|\psi|^2\big|^2}{|\psi|^4}
\end{align*}
which follows by a direct calculation using \eqref{eq:SchroedingerLichnerowiczmagnetic}.
Now, since \(\psi\in\Gamma(\Sigma^+M)\) is harmonic with respect to the magnetic Dirac operator \(D^\eta\) we can use 
the same strategy as in \cite[Proposition 3.2]{Branding:2018a} to establish that 
\(\big|d|\psi|^2\big|^2=2|\psi|^2|\nabla^\eta\psi|^2\).
Moreover, using the skew-symmetry of both \(d\eta\) and Clifford multiplication we find
\[
(d\eta)\cdot\psi=2e_1\cdot e_2\cdot (d\eta)(e_1,e_2)\psi
=-2d\eta(e_1,e_2)i\omega_\mathbb{C}\cdot \psi.
\]
By assumption we have that \(\psi\in\Gamma(\Sigma^\pm M)\) such that we are left with
\begin{align*}
-\Delta\log|\psi|=\frac{K}{2}\pm d\eta(e_1,e_2)
\end{align*}
and the claim follows by integrating over \(M\).
\end{proof}

\section{Diamagnetic inequality}\label{s:diamagneticinequality}
In this section, we give an obstruction for the diamagnetic inequality to hold. 
Recall that this inequality relates the eigenvalues of the magnetic Dirac operator to those with vanishing magnetic field.
We also investigate this inequality on Einstein Sasakian manifolds.
\begin{proposition}\label{prop:diamine} Let $(M^n,g)$ be a closed Riemannian spin manifold, $\eta\in \Omega^1(M)$ and $\psi$ 
be an eigenspinor of the Dirac operator associated with the eigenvalue $\lambda$. 
Then, for any $t\in\mathbb{R}$, we have 
$$(\lambda_1^{t\eta})^2\leq \lambda^2-t\frac{\int_M\Im\langle d\eta\cdot\psi-2\nabla_\eta\psi,\psi\rangle d\mu_g}{\|\psi\|^2_{L^2(M)}}+t^2\|\eta\|_\infty^2,$$
where $\Im$ resp. $\Re$ denote the imaginary resp. real parts of a complex number.
In particular, if $\int_M\Im\langle d\eta\cdot\psi-2\nabla_\eta\psi,\psi\rangle d\mu_g>0$, for some eigenspinor $\psi$, then
$$
|\lambda_1^{t\eta}|\leq |\lambda|,$$
holds for small positive parameter $t$.
\end{proposition}
\begin{proof}
Using the min-max principle, we compute for any eigenspinor field $\psi$ the following
\begin{eqnarray*}
\nonumber(\lambda_1^{t\eta})^2&\leq &\frac{\int_M\langle (D^{t\eta})^2\psi,\psi\rangle d\mu_g}{\|\psi\|^2_{L^2(M)}}\\
&\stackrel{\eqref{eq:relationdiracmagnwithdira}}{\leq} &\lambda^2+t\frac{\int_M \Re i\langle d\eta\cdot\psi-2\nabla_{\eta}\psi,\psi\rangle d\mu_g}{\|\psi\|^2_{L^2(M)}}+t^2\|\eta\|_\infty^2\\
\nonumber&=&\lambda^2-t\frac{\int_M\Im\langle d\eta\cdot\psi-2\nabla_\eta\psi,\psi\rangle d\mu_g}{\|\psi\|^2_{L^2(M)}}+t^2\|\eta\|_\infty^2.
\end{eqnarray*}
Hence, if $\int_M\Im\langle d\eta\cdot\psi-2\nabla_\eta\psi,\psi\rangle d\mu_g>0$, we deduce that, for small $t$, the r.h.s. of the above inequality is less than $\lambda^2$. 
\end{proof}

We now state the following

\begin{corollary}\label{cor:diaeinsasaki}
Let $(M^{2m+1},g,\eta)$ with $m>1$ be any simply connected closed $\eta$-Einstein Sasakian spin manifold of constant scalar curvature equal to $2m(2m-4b+1)$ with ($b<\frac{2m+1}{2}$) or ($b> \frac{2m+1}{2}$ and $m$ is odd). 
For any small positive $t$, we have 
$$|\lambda_1^{t\eta}|< |\frac{2m+1}{2}-b|.$$ 
In particular, when $b=0$, we deduce that $$|\lambda_1^{t\eta}|<\frac{2m+1}{2}=|\lambda_1|,$$ 
meaning that the diamagnetic inequality does not hold on simply connected closed Einstein Sasakian manifolds. 
\end{corollary} 

\begin{proof}
As mentioned at the beginning of Section \ref{sec:estimates}, any simply connected $\eta$-Einstein Sasakian manifold admits a Sasakian quasi-Killing spinor $\psi_m$ of type $(-\frac{1}{2},b)$ with $b\in \mathbb{R}$. The same holds for $\psi_0\in \Sigma_0M$ when $m$ is odd.
Hence $\psi_m$ (resp. $\psi_0$ when $m$ is odd) is an eigenspinor associated with the eigenvalue $\frac{2m+1}{2}-b$. 
Hence, we compute 
\begin{eqnarray*}
\int_M\Im\langle d\eta\cdot\psi_m-2\nabla_\eta\psi_m,\psi_m\rangle d\mu_g&=&\int_M\Im\langle 2im\psi_m-2(-\frac{1}{2}+b)\eta\cdot\psi_m,\psi_m\rangle d\mu_g\\
&=&(2m+1-2b)\int_M|\psi_m|^2 d\mu_g>0,
\end{eqnarray*}
when $b<\frac{2m+1}{2}$. 
In the last equality, we use the fact that $\eta\cdot\psi_m=i\psi_m$. 
When $m$ is odd, we have using a similar strategy
$$
\int_M\Im\langle d\eta\cdot\psi_0-2\nabla_\eta\psi_0,\psi_0\rangle d\mu_g=-(2m+1-2b)\int_M|\psi_0|^2d\mu_g.$$ 
Hence the condition in Proposition \ref{prop:diamine} is satisfied and, therefore we deduce, for small positive $t$, the required inequality. 
The last part comes from the fact that $\frac{2m+1}{2}$ is the lowest eigenvalue of the Dirac operator in absolute value. 
\end{proof}

For the $3$-dimensional case, we have the following 
\begin{corollary} 
Let $\mathbb{S}^3$ be the $3$-dimensional sphere with metric of the form $g:=s^2 \eta\otimes\eta+s g_0{}_{|_{\eta^\perp}}$ for some $s\in(0,\infty)$, where $\eta$ is the Reeb one-form and $g_0$ is the standard round metric of constant sectional curvature $1$.
Denote by $S$ its (constant) scalar curvature.
Then, for any sufficiently small positive parameter $t$, the following estimate holds
$$|\lambda_1^{t\eta}|<|\frac{3}{4}+\frac{S}{8}|.$$ 
\end{corollary} 
\begin{proof}
As before, it is shown in \cite[Thm. 8.4]{FriedrichKim:00} that there exist two Sasakian quasi-Killing spinors $\psi_0$ and $\psi_1$ on $\Sigma_0 M$ and $\Sigma_1 M$ of type $(-\frac{1}{2},\frac{3}{4}-\frac{1}{8}S)$.
Hence, both $\psi_0$ and $\psi_1$ are eigenspinors of the Dirac operator for the eigenvalue $\frac{3}{4}+\frac{1}{8}S$. 
As before, we have  
\begin{eqnarray*}
\int_M\Im\langle d\eta\cdot\psi_1-2\nabla_\eta\psi_1,\psi_1\rangle d\mu_g&=&(\frac{3}{2}+\frac{1}{4}S)\int_M|\psi_1|^2 d\mu_g,
\end{eqnarray*}
and $\int_M\Im\langle d\eta\cdot\psi_0-2\nabla_\eta\psi_0,\psi_0\rangle d\mu_g=-(\frac{3}{2}+\frac{1}{4}S)\int_M|\psi_0|^2 d\mu_g$. 
This concludes the proof.
\end{proof}

\section{Computations of the spectrum}\label{s:computmagneticspectra}
In this section, we compute explicitly the spectrum of the magnetic Dirac operator on the round $3$-sphere and on the $n$-dimensional flat torus. 

\subsection{The 3-dimensional round sphere}\label{ss:computspectrumS3}
The spectrum of the magnetic Dirac operator on the 3-dimensional round sphere associated with the Reeb vector field can be described as follows, see also \cite{ErdoesSolovej:01} for a related computation.

\begin{theorem}\label{t:spectrumDetaS3}
Let $\mathbb{S}^3$ be the round sphere equipped with the metric of curvature $1$. 
Let $\eta$ be the Reeb vector field that defines the Hopf fibration. 
Then, for any $t\in \mathbb{R}$, the spectrum of the magnetic Dirac operator $D^{t\eta}$ is given by
$$\frac{3}{2}\pm t+k, \frac{1}{2}\pm \sqrt{(1+t+2p-k)^2+4(k-p)(p+1)},$$
for $p,k\in \NN$ with $0\leq p<k$.
\end{theorem}
\begin{proof} Recall that the Hopf vector field on the round sphere $\mathbb{S}^3\subset\mathbb{C}^2$ equipped with the metric $g$ of curvature $1$ is given by $\eta=\sum_{j=1}^2(-y_j\partial_{x_j}+x_j\partial_{y_j})$. Next, we follow the computations done in \cite{EGHP:23}  for the spectrum of the magnetic Laplacian (see also \cite{norbert-thesis}). We let $Y_2=\eta, Y_3, Y_4$ to be the Killing vector fields  given by
\begin{eqnarray*}
Y_2 &=& -y_1 \partial_{x_1} + x_1 \partial_{y_1} - y_2 \partial_{x_2} + x_2 \partial_{y_2}, \\
Y_3 &=& -y_2 \partial_{x_1} - x_2 \partial_{y_1} + y_1 \partial_{x_2} + x_1 \partial_{y_2}, \\
Y_4 &=& x_2 \partial_{x_1} - y_2 \partial_{y_1} - x_1 \partial_{x_2} + y_1 \partial_{y_2}.
\end{eqnarray*}
They form a direct orthonormal basis of $T_{(z_1,z_2)}\mathbb{S}^3$ at every point $(z_1,z_2) = (x_1 + y_1 i, x_2 + y_2 i) \in \mathbb{S}^3$. An easy computation shows that the Christoffel symbols of the Levi-Civita connection of $g$ are expressed as
\begin{equation}\label{eq:chrsym} 
\nabla_{Y_j} Y_k= \sigma_{jk} Y_l
\end{equation}
with $\{j,k,l\} = \{2,3,4\}$ for $k \neq j$, $\sigma_{jj} = 0$ and $\sigma_{23}=-\sigma_{24} =-1$, $\sigma_{32}=\sigma_{43} = - \sigma_{34} = -\sigma_{42} = 1$. 
Recall now that the eigenvalues of the scalar Laplacian are given by $k(k+2)$ for $k\in\NN\cup \{0\}$ with multiplicity $(k+1)^2$. 
Each eigenspace $E_k$ can be decomposed as 
\begin{equation} \label{eq:Ekdecomp}
E_k= V_{k, (a_0,b_0)} \oplus V_{k, (a_1, b_1)}\oplus \ldots \oplus V_{k, (a_k, b_k)},
\end{equation}
with any arbitrary choice of pairwise non-collinear vectors $(a_j,b_j) \in \CC \setminus\{(0,0)\}$, where
\begin{align*}
& V_{k, (a,b)}  = {\textrm{span}}_\CC\{u_{a,b}^k, u_{a,b}^{k-1}v_{a,b}, \ldots, u_{a,b} v_{a,b}^{k-1}, v_{a,b}^k\}, \\
& u_{a,b}(z_1, z_2):= az_1+bz_2, \quad v_{a,b}(z_1, z_2):=b \bar z_1 -a\bar z_2,
\end{align*}
for $(a,b)\in\CC^2\setminus\{(0,0)\}$, see \cite[Zerlegungssatz III.6.2]{norbert-thesis}. For short, we write $u:= u_{(a,b)}$, $v:= v_{(a,b)}$ for some $(a,b)\neq (0,0)$ and, for $p \in \{0,\ldots, k\}$, we consider
$$
\phi_{k,p}:= u^p v^{k-p}.
$$
We also set $\phi_p \equiv 0$ for all other choices of $p$. These functions $\phi_{k,p}$ are spherical harmonics, that is, they are restrictions of harmonic homogeneous polynomials on $\CC^2$ to the unit sphere $\mathbb{S}^3$. Then we have $V_{k, (a,b)}={\textrm{span}}_\CC\{\phi_0,\ldots, \phi_k\}$.
A straightforward computation yields (see \cite[p. 30]{Hitchin:74} or \cite[Lemma III.7.1]{norbert-thesis})
\begin{eqnarray}\label{eq:Y2phi}
Y_2(\phi_{k,p}) &=& i(2p-k)\phi_{k,p},\nonumber\\
Y_3(\phi_{k,p}) &=& ip\phi_{k,p-1}+i(k-p)\phi_{k,p+1}, \nonumber\\
Y_4(\phi_{k,p}) &=& -p\phi_{k,p-1}+(k-p)\phi_{k,p+1}.\nonumber 
\end{eqnarray}
It is shown in \cite[Theorem 8.4]{FriedrichKim:00} that the spinor bundle of $\mathbb{S}^3$ is trivialised by $-\frac{1}{2}$-Killing spinors which are pointwise eigenvectors for the Clifford action of $\eta$, all terms -- including that involving the Clifford action of the differential of the function $\phi_{k,p}$ -- can be explicitly computed.
Let $\psi_0,\psi_1$ be $-\frac{1}{2}$-Killing spinors of unit length on $\mathbb{S}^3$ with $i\eta\cdot\psi_j=(-1)^j\psi_j$, $j=0,1$. We can assume that $Y_3\cdot\psi_0=\psi_1$. Indeed, the spinor $Y_3\cdot\psi_0$ is in $\Sigma_1M$ which is also a $-\frac{1}{2}$-Killing spinor. To see this, we compute 
\begin{eqnarray*}
\nabla_\eta(Y_3\cdot\psi_0)&=&\nabla_\eta Y_3\cdot\psi_0+Y_3\cdot\nabla_{\eta}\psi_0\\
&\stackrel{\eqref{eq:chrsym}}{=}&-Y_4\cdot\psi_0-\frac{1}{2}Y_3\cdot\eta\cdot\psi_0\\
&=&-iY_3\cdot\psi_0-\frac{1}{2}Y_3\cdot\eta\cdot\psi_0\\
&=&\frac{1}{2}Y_3\cdot\eta\cdot\psi_0=-\frac{1}{2}\eta\cdot Y_3\cdot\psi_0.
\end{eqnarray*}
In this computation, we use the fact that $Y_4\cdot\psi_0=iY_3\cdot\psi_0$, since $Y_4=JY_3$, where $J$ is the complex structure on $\eta^\perp$. In the same way, we compute 
\begin{eqnarray*}
\nabla_{Y_3}(Y_3\cdot\psi_0)&=&\nabla_{Y_3} Y_3\cdot\psi_0+Y_3\cdot\nabla_{Y_3}\psi_0\\
&\stackrel{\eqref{eq:chrsym}}{=}&-\frac{1}{2}Y_3\cdot Y_3\cdot\psi_0, 
\end{eqnarray*}
and 
\begin{eqnarray*}
\nabla_{Y_4}(Y_3\cdot\psi_0)&=&\nabla_{Y_4} Y_3\cdot\psi_0+Y_3\cdot\nabla_{Y_4}\psi_0\\
&\stackrel{\eqref{eq:chrsym}}{=}&\eta\cdot\psi_0-\frac{1}{2}Y_3\cdot Y_4\cdot\psi_0\\
&=&-\frac{1}{2}Y_4\cdot Y_3\cdot\psi_0.
\end{eqnarray*}
In the last equality, we use the fact that $\eta\cdot\psi_0=Y_3\cdot Y_4\cdot\psi_0$ since the volume form $-\eta\cdot Y_3\cdot Y_4\cdot$ acts by the identity on $\Sigma \mathbb{S}^3$. 
As the space of $-\frac{1}{2}$-Killing spinors in $\Sigma_1 M$ is one-dimensional, we deduce that $Y_3\cdot\psi_0$ is collinear to $\psi_1$. Hence up to rescaling $\psi_0$, we can assume that $Y_3\cdot\psi_0=\psi_1$. For $j\in\{0,1\}$, we have that $D^{t\eta}\psi_j=D\psi_j+it\eta\cdot\psi_j=(\frac{3}{2}+(-1)^jt)\psi_j$ and, therefore, we compute
\begin{eqnarray*}
(D^{t\eta}-\frac{1}{2})(\phi_{k,p}\psi_j)&=&\phi_{k,p}(D^{t\eta}-\frac{1}{2})\psi_j+d\phi_{k,p}\cdot\psi_j\\
&=&\phi_{k,p}\left(\frac{3}{2}+(-1)^jt-\frac{1}{2}\right)\psi_j+d\phi_{k,p}\cdot\psi_j\\
&=&(1+(-1)^jt)\phi_{k,p}\psi_j+d\phi_{k,p}\cdot\psi_j.
\end{eqnarray*}
Now, we explicitly compute $d\phi_{k,p}\cdot\psi_j$.
From \eqref{eq:Y2phi}, we have that
\[d\phi_{k,p}=i(2p-k)\phi_{k,p}\eta+i(p\phi_{k,p-1}+(k-p)\phi_{k,p+1})Y_3+((k-p)\phi_{k,p+1}-p\phi_{k,p-1})Y_4,\]
and for both $j=0,1$ we find
\begin{eqnarray*}
d\phi_{k,p}\cdot\psi_j&=&i(2p-k)\phi_{k,p}\eta\cdot\psi_j+i(p\phi_{k,p-1}+(k-p)\phi_{k,p+1})Y_3\cdot\psi_j\\
&&+((k-p)\phi_{k,p+1}-p\phi_{k,p-1})Y_4\cdot\psi_j\\
&=&(-1)^j(2p-k)\phi_{k,p}\psi_j+(-1)^ji(p\phi_{k,p-1}+(k-p)\phi_{k,p+1})\psi_{j+1}\\
&&+i((k-p)\phi_{k,p+1}-p\phi_{k,p-1})\psi_{j+1}\\
&=&(-1)^j(2p-k)\phi_{k,p}\psi_j\\
&&+i\left((1+(-1)^j)(k-p)\phi_{k,p+1}+((-1)^j-1)p\phi_{k,p-1}\right)\psi_{j+1},
\end{eqnarray*}
where the index $j+1$ is to be understood mod $2$.
Note that we have used the fact that $Y_4\cdot\psi_0=iY_3\cdot\psi_0=i\psi_1$. Thus, it holds that
\begin{eqnarray*}
d\phi_{k,p}\cdot\psi_0&=&(2p-k)\phi_{k,p}\psi_0+2i(k-p)\phi_{k,p+1}\psi_1,\\
d\phi_{k,p}\cdot\psi_1&=&-(2p-k)\phi_{k,p}\psi_1-2ip\phi_{k,p-1}\psi_0.
\end{eqnarray*}
Therefore, we get
\begin{eqnarray*}
(D^{t\eta}-\frac{1}{2})(\phi_{k,p}\psi_0)&=&(1+t)\phi_{k,p}\psi_0+(2p-k)\phi_{k,p}\psi_0+2i(k-p)\phi_{k,p+1}\psi_1\\
&=&(1+t+2p-k)\phi_{k,p}\psi_0+2i(k-p)\phi_{k,p+1}\psi_1,
\end{eqnarray*}
and similarly
\[(D^{t\eta}-\frac{1}{2})(\phi_{k,p}\psi_1)=(1-t-2p+k)\phi_{k,p}\psi_1-2ip\phi_{k,p-1}\psi_0.\]
Note that, even if neither $\phi_{k,k+1}$ nor $\phi_{k,-1}$ are defined, we still have $\displaystyle(D^{t\eta}-\frac{1}{2})(\phi_{k,k}\psi_0)=(1+t+k)\phi_{k,k}\psi_0$ when $p=k$ and $\displaystyle(D^{t\eta}-\frac{1}{2})(\phi_{k,0}\psi_1)=(1-t+k)\phi_{k,0}\psi_1$ when $p=0$.
Note also that, when $p<k$ resp. $p>0$, the spinor field $\phi_{k,p}\psi_0$ resp. $\phi_{k,p}\psi_1$ is not an eigenvector for $D^{t\eta}-\frac{1}{2}$.\\

Now, for $0\leq p<k$, let $\varphi_{k,p,0}^\pm:=(\mathrm{Id}\pm\frac{1}{\sqrt{f_0(k,p,t)}}(D^{t\eta}-\frac{1}{2}))(\phi_{k,p}\psi_0)$, where $f_0(k,p,t):=(1+t+2p-k)^2+4(k-p)(p+1)\geq0$.
Because of 
\begin{eqnarray*}
(D^{t\eta}-\frac{1}{2})^2(\phi_{k,p}\psi_0)&=&(D^{t\eta}-\frac{1}{2})\left((1+t+2p-k)\phi_{k,p}\psi_0+2i(k-p)\phi_{k,p+1}\psi_1\right)\\
&=&(1+t+2p-k)(D^{t\eta}-\frac{1}{2})(\phi_{k,p}\psi_0)\\
&&+2i(k-p)(D^{t\eta}-\frac{1}{2})(\phi_{k,p+1}\psi_1)\\
&=&(1+t+2p-k)\left((1+t+2p-k)\phi_{k,p}\psi_0+2i(k-p)\phi_{k,p+1}\psi_1\right)\\
&&+2i(k-p)\left((1-t-2(p+1)+k)\phi_{k,p+1}\psi_1-2i(p+1)\phi_{k,p}\psi_0\right)\\
&=&\left((1+t+2p-k)^2+4(k-p)(p+1)\right)\phi_{k,p}\psi_0\\
&=&f_0(k,p,t)\phi_{k,p}\psi_0,
\end{eqnarray*}
the spinor fields $\varphi_{k,p,0}^\pm\neq0$ satisfy $(D^{t\eta}-\frac{1}{2})\varphi_{k,p,0}^\pm=\pm\frac{1}{\sqrt{f_0(k,p,t)}}\varphi_{k,p,0}^\pm$.
Similarly, one could define $\varphi_{k,p,1}^\pm:=(\mathrm{Id}\pm\frac{1}{\sqrt{f_1(k,p,t)}}(D^{t\eta}-\frac{1}{2}))(\phi_{k,p}\psi_1)$ for all $0<p\leq k$, where $f_1(k,p,t):=(1-t-2p+k)^2+4p(k-p+1)$.
A fundamental remark is that $\varphi_{k,p,1}^\pm$ is a scalar multiple of $\varphi_{k,p-1,0}^\pm$: namely $f_1(k,p,t)=f_0(k,p-1,t)$ by definition and the matrix expressing the pair $(\varphi_{k,p,0}^\pm,\varphi_{k,p+1,1}^\pm)$ in terms of $(\phi_{k,p}\psi_0,\phi_{k,p+1}\psi_1)$ has rank at most $1$. On the whole, we have found three families of eigenvectors and eigenvalues for $D^{t\eta}-\frac{1}{2}$: 
$$\{\phi_{k,k}\psi_0\,|\,k\in\mathbb{N}\},\,\,\,\,\,  \{\phi_{k,0}\psi_1\,|\,k\in\mathbb{N}\} \quad\text{and}\quad \{\varphi_{k,p,0}^\pm\,|\,p,k\in\mathbb{N}, 0\leq p<k\}$$ 
with the respective eigenvalues being $1+t+k$, $1-t+k$ and $\pm\sqrt{f_0(k,p,t)}$.
Note that $(1+t+k)^2=f_0(k,k,t)$ and $(1-t+k)^2=f_0(k,-1,t)$, so that all eigenvalues actually have the same expression.
The union of those three families forms an orthogonal Hilbert basis of $L^2(\Sigma M)$ since the family $\{\phi_{k,p}\psi_0,\phi_{k,p}\psi_1\}$ does. This shows the spectrum of $D^{t\eta}$ as described in the theorem.\\

It remains to determine the multiplicity of each eigenvalue.
Since, for each $p,k\in\mathbb{N}$ with $0\leq p\leq k$, the functions $\phi_{k,p}$ -- which depend on $u$ and $v$ and therefore also on $(a,b)\in\mathbb{C}^2\setminus\{(0,0)\}$ -- form a $(k+1)$-dimensional space, as we have seen above, each of the above eigenvalue has multiplicity at least $k+1$.
Beware here that $\varphi_{k,p,0}^\pm$ is uniquely determined by $\phi_{k,p}\psi_0$ and that the eigenvalues of two of the three families may coincide for particular values of $t$.
For instance, $f_0(k,p,t)=f_0(k',p',t)$ for $k,p,k',p'$ as above if and only if $(1+t+2p-k)^2+4(k-p)(p+1)=(1+t+2p'-k')^2+4(k'-p')(p'+1)$, which is equivalent to $(t+p+p'-\frac{k+k'}{2}+1)(p-p'-\frac{k-k'}{2})+(k-p)(p+1)-(k'-p')(p'+1)=0$.
If $p-p'=\frac{k-k'}{2}$ i.e., $k'=k+2(p'-p)$ then $(k-p)(p+1)-(k'-p')(p'+1)=0$ is equivalent to $(k-p)(p+1)-(k+p'-2p)(p'+1)=0$, which amounts to $(p-p')(k+1-p+p')=0$, from which either $p=p'$ and then $k=k'$ follows, or $k=p-p'-1$, which is a contradiction because of $p-p'-1<p<k$.
Therefore, if $(k,p)\neq(k',p')$, then $p-p'\neq\frac{k-k'}{2}$ and there exists a unique solution $t$ to $f_0(k,p,t)=f_0(k',p',t)$ which is given by \[t=\frac{(k'-p')(p'+1)-(k-p)(p+1)}{p-p'-\frac{k-k'}{2}}+\frac{k+k'}{2}-p-p'-1.\]
The other cases are analogous.
\begin{figure}[ht]
    \centering
    \begin{tikzpicture}
\clip (-5,-5) rectangle (5,5);
\draw[step=.5cm,gray,very thin] (-5,-5) grid (5,5);
\draw (0,0) node[fill=white,anchor=north east] {$0$};
\foreach \x in {-4,-3,-2,-1,1,2,3,4}
\draw (\x,0) node[fill=white,anchor=north] {$\x$};
\foreach \x in {-4,-3,-2,-1,1,2,3,4}
\draw[thick] (\x,-0.1) -- (\x,0.1);
\foreach \x in {-4,-3,-2,-1,1,2,3,4}
\draw (0,\x) node[fill=white,anchor=east] {$\x$};
\foreach \x in {-4,-3,-2,-1,1,2,3,4}
\draw[thick] (-0.1,\x) -- (0.1,\x);
\foreach \k in {0,...,5}
\draw[thick,color=blue,smooth,samples=100,domain=-5:5] plot (\x,{1.5+\k+\x});
\foreach \k in {0,...,5}
\draw[thick,color=red,smooth,samples=100,domain=-5:5] plot (\x,{1.5+\k-\x});
\foreach \k in {1,...,4}
 \foreach[parse=true] \p in {0,...,\k-1}
\draw[thick,color=green!60!black,smooth,samples=100,domain=-5:5] plot (\x,{0.5+sqrt((\x+1+2*\p-\k)^2+4*(\k-\p)*(\p+1))});
\foreach \k in {1,...,4}
 \foreach[parse=true] \p in {0,...,\k-1}
 \draw[thick,color=green!60!black,smooth,samples=100,domain=-5:5] plot (\x,{0.5-sqrt((\x+1+2*\p-\k)^2+4*(\k-\p)*(\p+1))});
\draw (4.8,0.05) node[fill=white,anchor=south] {$t$};
\draw[->,thick] (-5,0) -- (5,0) coordinate (x axis);
\draw[->,thick] (0,-5) -- (0,5) coordinate (y axis);
\end{tikzpicture}
    \caption{Eigenvalues of $D^{t\eta}$ on $\mathbb{S}^3$ as functions of $t$: {\color{blue}$\frac{3}{2}+k+t$}, {\color{red}$\frac{3}{2}+k-t$}, {\color{green!60!black}$\frac{1}{2}\pm\sqrt{f_0(k,p,t)}$}}
    \label{fig:eigenvaluesofS3}
\end{figure}
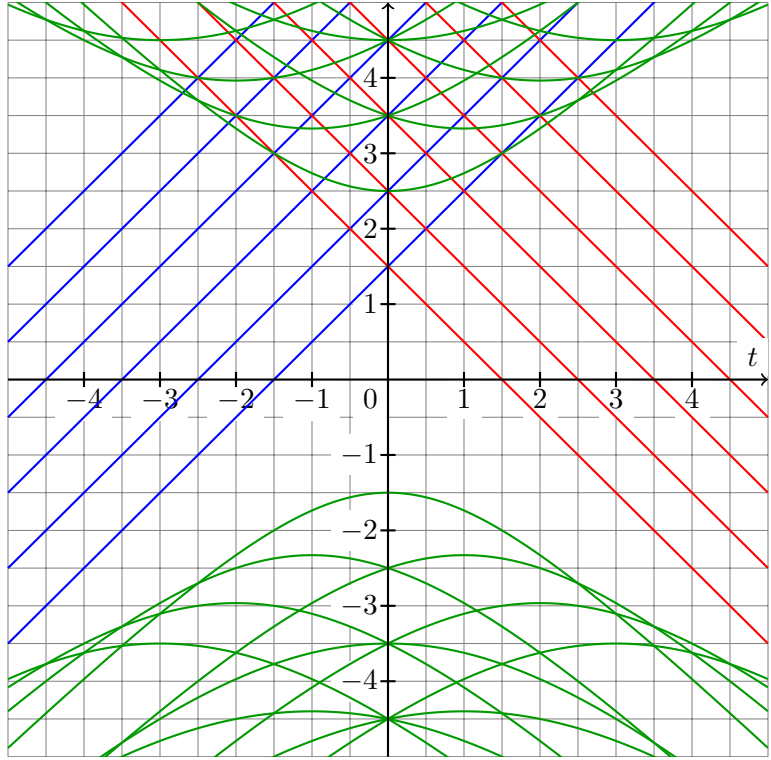
For $t=0$, the multiplicity of the eigenvalue $1+k$ of $D^{t\eta}-\frac{1}{2}$ is $(k+2)(k+1)$ as expected \cite[Sec. 2]{Baerpsf:96} since there is a $k+1$-dimensional contribution from each of the eigenvectors $\phi_{k,k}\psi_0,\phi_{k,0}\psi_1,\varphi_{k,p,0}^+$, $0\leq p<k$.\\
\end{proof}
Let us make some comments concerning the eigenvalue estimates discussed in the previous sections.
\begin{remark}\label{rem:magneticspectrumS3}\rm
\noindent\begin{enumerate} 
\item The magnetic Friedrich inequality in Theorem \ref{thm:magneticfrie} cannot be optimal on the round sphere, since the latter is an Einstein Sasakian manifold. However, we can check this by the computation of the spectrum done for $\mathbb{S}^3$. Namely, for $t\in ]0,\frac{3}{2}]$, we have $|\lambda_1^{t\eta}|=\frac{3}{2}-t$ and 
$$\frac{n}{4(n-1)}\mathop{\rm inf}\limits_M(S-4t\lfloor\frac{n}{2}\rfloor^\frac{1}{2}|d\eta|)=\frac{3}{4}(3-4t)$$
which is clearly strictly less than $(\frac{3}{2}-t)^2$.
\item The inequality $|\lambda_1^{t\eta}|<|\lambda_1|=\frac{3}{2}$ from  Corollary \ref{cor:diaeinsasaki} holds \emph{for any} $t\in\mathbb{R}\setminus\{0\}$ and not only for small $t$. Indeed, the eigenvalue $\frac{3}{2}-t$ of $D^{t\eta}$ satisfies that inequality for all $t\in(0,3)$.
But so does $\frac{5}{2}-t$ for all $t\in(1,4)$.
Inductively, it can be shown that the $D^{t\eta}$-eigenvalue $\frac{3}{2}+k-t$ satisfies that inequality for all $t\in(k,k+3)$. Therefore, $|\lambda_1^{t\eta}|<|\lambda_1|=\frac{3}{2}$ holds for all positive $t$. Replacing $t$ by $-t$ and using the $D^{t\eta}$-eigenvalues $\frac{3}{2}+k+t$, the same holds for all negative $t$, hence the conclusion.
\item The inequality \eqref{eq:Hijazimagneticineq} is an equality for $\lambda_1^{t\eta}$ if and only if $t\in[-\frac{3}{2},\frac{3}{2}]$.
Note that $\|t\eta\|_{L^3(\mathbb{S}^3)}=|t|\omega_3^{\frac{1}{3}}$, where $\omega_3$ denotes the volume of the round $3$-sphere $(\mathbb{S}^3,g)$.
On the other hand, $\sqrt{\frac{3}{4(3-1)}Y(\mathbb{S}^3,[g])}=\sqrt{\frac{3}{8}\cdot3\cdot2\cdot\omega_3^{\frac{2}{3}}}=\frac{3}{2}\omega_3^{\frac{1}{3}}$.
Therefore, \eqref{eq:Hijazimagneticineq} is an equality if and only if $|\lambda_1^{t\eta}|+|t|=\frac{3}{2}$.
For $t\in[0,\frac{3}{2}]$, $|\lambda_1^{t\eta}|=\frac{3}{2}-t$, so that $|\lambda_1^{t\eta}|+|t|=\frac{3}{2}$.
Analogously, for $t\in[-\frac{3}{2},0]$, $|\lambda_1^{t\eta}|=\frac{3}{2}+t$, so that $|\lambda_1^{t\eta}|+|t|=\frac{3}{2}$.
For $|t|>\frac{3}{2}$, obviously $|\lambda_1^{t\eta}|+|t|>\frac{3}{2}$ holds and therefore \eqref{eq:Hijazimagneticineq} cannot be an equality.
 Note that $|\lambda^{t\eta}|\leq\frac{1}{2}$ for all $|t|\geq\frac{3}{2}$, see e.g. figure \ref{fig:eigenvaluesofS3}.
\item When $M:=\rquot{\mathbb{S}^3}{\Gamma}$ with round metric $g$ of sectional curvature $1$, where $\Gamma$ is any finite subgroup of $\mathrm{SU}_2$, 
there exists a spin structure for which the space of $\frac{1}{2}$-Killing spinors is complex $2$-dimensional and that of $-\frac{1}{2}$-Killing spinors $0$-dimensional; and a spin structure for which the same holds true where the constants $\pm\frac{1}{2}$ are swapped \cite[Cor. 5.2.5]{Ammann:habil}.
Since, as mentioned above, the spinor bundle of $\mathbb{S}^3$ is trivialised by $-\frac{1}{2}$-Killing spinors which are pointwise $\pm1$-eigenvectors of $i\eta\cdot$ \cite[Theorem 8.4]{FriedrichKim:00}, there exists a spin structure on $M$ for which the spinor bundle $\Sigma M$ is trivialised by $-\frac{1}{2}$-Killing spinors on $M$ which are pointwise $\pm1$-eigenvectors of $i\eta\cdot$.
Note that $\eta$ is anyway $\Gamma$-invariant since $\Gamma$ acts on $\mathbb{C}^2$ by \emph{unitary} (thus complex-linear) isomorphisms.
That spin structure being fixed, we obtain $D^{t\eta}\psi=0$ for any such $-\frac{1}{2}$-Killing spinor $\psi$ and for $t\in\{\pm\frac{3}{2}\}$.
Therefore, after using $\mathrm{Vol}(M,g)=\frac{\omega_3}{|\Gamma|}$, \eqref{eq:Hijazimagneticineq} becomes 
\begin{equation}\label{eq:upperboundYS3overGamma}
Y(\rquot{\mathbb{S}^3}{\Gamma},[g])\leq\frac{Y(\mathbb{S}^3,[g])}{|\Gamma|^{\frac{2}{3}}}.
\end{equation}
When $\Gamma=\{\pm1\}$ i.e., when $M$ is the $3$-dimensional real projective space $\mathbb{R}\mathrm{P}^3$, it is known though nontrivial that \eqref{eq:upperboundYS3overGamma} is an equality \cite{BrayNeves:2004}.
Therefore, \eqref{eq:Hijazimagneticineq} is also an equality on $M=\mathbb{R}\mathrm{P}^3$, in particular $\mathbb{S}^3$ is not the only closed Riemannian $3$-manifold on which \eqref{eq:Hijazimagneticineq} is an equality, in spite of the claim of \cite[Remark 25]{Reuss:25}.
For further groups $\Gamma$, it is still an open question whether \eqref{eq:upperboundYS3overGamma} is an equality or not.
\end{enumerate}
\end{remark}
We end this section by pointing out that the claim of \cite[Theorem 1.1]{FrankLoss:2021} in dimension $3$ -- as well as of \cite[Theorem 1.3]{FrankLoss:2021} in higher dimensions -- is likely to be implied by \cite[Lemma 3.4]{HerzlichMoroianu:99} using a stereographic projection to map $\RR^n$ into $\mathbb{S}^n$ and the conformal covariance of the Dirac operator as above.

\subsection{The flat torus}\label{ss:specspincflattorus}
In this section, we compute the spectrum of the spin$^c$ Dirac operator on the flat torus $\mathbb{T}^n$,
the case of the flat three torus was already studied in \cite{Meier:12} and the spectrum of the spin Dirac operator on flat $\mathbb{T}^n$ was computed in \cite{Friedrich:84}.
We start with the case where the auxiliary one-form is parallel, the case where that one-form is only assumed to be closed being discussed after the proof of Theorem \ref{t:specflattorusspincAparallel}.
We emphasize that, because every spin$^c$ structure on the torus actually reduces to a spin structure, spin$^c$ Dirac operators on $\mathbb{T}^n$ coincide with magnetic Dirac operators.
The notations needed to understand the statement of Theorem \ref{t:specflattorusspincAparallel} will be introduced in its proof.

\begin{theorem}\label{t:specflattorusspincAparallel}
For a positive integer $n$ let $\Gamma\subset\RR^n$ be a lattice and $\mathbb{T}^n:=\lquot{\RR^n}{\Gamma}$ the corresponding $n$-dimensional flat torus.
Let $P_{\mathbb{U}_1}\to M$ be a fixed $\mathbb{U}_1$-bundle and let a square root of that bundle be given.
Fix a basis $(\gamma_1,\ldots,\gamma_n)$ of $\Gamma$ and $\delta_1,\ldots,\delta_n\in\{0,1\}$.
Denote by $\theta_1,\ldots,\theta_n\in[0,1[$ the real numbers fixing the square root of $P_{\mathbb{U}_1}\to M$.
Let $A$ be a parallel one-form on $\mathbb{T}^n$. 
Then the spectrum of the Dirac operator of $\mathbb{T}^n$ endowed with the induced flat metric and the spin$^c$-structure provided by $(\delta_1,\ldots,\delta_n)$ and $(\theta_1,\ldots,\theta_n)$ is given by
\[ \Big\{\pm 2\pi|\gamma^*+\frac{1}{2}\sum_{j=1}^n(\delta_j+\theta_j)\gamma_j^*+\frac{A}{4\pi}|,\quad\gamma^*\in\Gamma^*\Big\},\]
where $\Gamma^*:=\{\theta\in(\RR^n)^*\,|\,\theta(\Gamma)\subset\mathbb{Z}\}$ is the dual lattice and $(\gamma_1^*,\ldots,\gamma_n^*)$ the basis of $\Gamma^*$ dual to $(\gamma_1,\ldots,\gamma_n)$. 
If non-zero, the eigenvalue provided by $\gamma^*$ has multiplicity at least $2^{[\frac{n}{2}]-1}$.
Furthermore, $0$ is an eigenvalue if and only if there exists $\gamma^*\in\Gamma^*$ such that $A=-4\pi(\gamma^*+\frac{1}{2}\sum_{j=1}^n(\delta_j+\theta_j)\gamma_j^*)$ and, in that case, its multiplicity is exactly $2^{[\frac{n}{2}]}$ and the corresponding eigenspace consists of $\nabla^A$-parallel spinors.
\end{theorem}
\begin{proof} Let $P_{\mathbb{U}_1}\to \mathbb{T}^n$ be any auxiliary $\mathbb{U}_1$-bundle on $\mathbb{T}^n$, which can be described by a group-homomorphism $\beta\colon\Gamma\to\mathbb{U}_1$.
A spin$^c$ structure on $\mathbb{T}^n$ with associated $\mathbb{U}_1$-bundle provided by $\beta$ is then fully described by a group-homomorphism $\hat{\beta}\colon\Gamma\to\mathrm{Spin}_n^c$ such that, for every $\gamma\in\Gamma$, we have $\eta^c(\hat{\beta}(\gamma))=(1,\beta(\gamma))$.
Here $\eta^c\colon\mathrm{Spin}_n^c\to\mathrm{SO}_n\times\mathbb{U}_1$, is the nontrivial $2$-fold covering.
The reason for the first component of the group-homomorphism $(1,\beta)$ to be $1$ is that $\Gamma$ acts trivially on the tangent bundle of $\mathbb{R}^n$, since it consists of translations only.
Now recall that $\mathrm{Spin}_n^c\cong\rquot{\mathrm{Spin}_n\times\mathbb{U}_1}{\mathbb{Z}_2}$ via the surjective map $\mathrm{Spin}_n\times\mathbb{U}_1\to\mathrm{Spin}_n^c$, $(u,z)\mapsto zu$ with kernel $\{\pm(1,1)\}$; and that, via that identification, $\eta^c([(u,z)])=(\eta(u),z^2)$ for all $(u,z)\in\mathrm{Spin}_n\times\mathbb{U}_1$, where $\eta\colon\mathrm{Spin}_n\to\mathrm{SO}_n$ is the nontrivial twofold covering.
Fix $\gamma\in\Gamma$ and write $\hat{\beta}(\gamma)=[(u(\gamma),z(\gamma))]$ for some $u(\gamma)\in\mathrm{Spin}_n$ and $z(\gamma)\in\mathbb{U}_1$.
Then the identity $\eta^c(\hat{\beta}(\gamma))=(1,\beta(\gamma))$ is equivalent to $\eta(u(\gamma))=1$ and $z(\gamma)^2=\beta(\gamma)$.
But because of $\ker(\eta)=\{\pm1\}$, the identity $\eta(u(\gamma))=1$ is equivalent to $u(\gamma)\in\{\pm 1\}\subset\mathrm{Spin}_n$.
Therefore, $\hat{\beta}$ can be written in two ways.
Either a square-root $\tilde{\beta}\colon\Gamma\to\mathbb{U}_1$ of $\beta$ is fixed and then $\hat{\beta}(\gamma)=[(\varepsilon(\gamma),\tilde{\beta}(\gamma))]$ for all $\gamma\in\Gamma$ and some group-homomorphism $\varepsilon\colon\Gamma\to\{\pm1\}$; or the first component of $\hat{\beta}$ is fixed to $1$ i.e., $\hat{\beta}(\gamma)=[(1,\tilde{\beta}'(\gamma))]$ for all $\gamma\in\Gamma$ and some group-homomorphism $\tilde{\beta}'\colon\Gamma\to\mathbb{U}_1$ which is a square-root of $\beta$.
One way or the other, there are as many spin$^c$ structures on $\mathbb{T}^n$ with auxiliary $\mathbb{U}_1$-bundle given by $\beta$ as group-homomorphisms $\Gamma\to\{\pm1\}$.
In particular, there are $2^n$ such spin$^c$ structures on $\mathbb{T}^n$ since $\Gamma$ has a basis consisting of $n$ linearly independent vectors of $\RR^n$.
Note that the existence of a square-root of $\beta$ is anyway ensured by $\Gamma$ being finitely generated (define the square-root on a chosen basis and extend it as a group-homomorphism on the whole $\Gamma$).\\

As a consequence, any spin$^c$ structure on $\mathbb{T}^n$ reduces to a spin structure via the above map $\mathrm{Spin}_n\times\mathbb{U}_1\to\mathrm{Spin}_n^c$.
In the former description where the square-root $\tilde{\beta}$ of $\beta$ is fixed, the spin$^c$ structures on $\mathbb{T}^n$ stand in one-to-one correspondence with spin structures on $\mathbb{T}^n$, a square-root-bundle of the auxiliary bundle being fixed.
In the latter description where the $\mathrm{Spin}_n$-component of $\hat{\beta}$ is fixed to $1$, the spin$^c$ structures on $\mathbb{T}^n$ stand in one-to-one correspondence to the trivial spin structure cross a square-root-bundle of the auxiliary bundle.\\

Let the square-root $\tilde{\beta}\colon\Gamma\to\mathbb{U}_1$ of $\beta$ be fixed and let $\varepsilon\colon\Gamma\to\{\pm1\}$ describe the spin$^c$ structure on $\mathbb{T}^n$.
Then a spinor field $\phi$ on $\mathbb{T}^n$ may be identified with a spinor field $\phi$ on $\RR^n$ satisfying, in the above notations,
\[\phi(x+\gamma)=\varepsilon(\gamma)\tilde{\beta}(\gamma)\phi(x),\] 
for all $x\in\RR^n$ and $\gamma\in\Gamma$.\\

Let $(\gamma_1,\ldots,\gamma_n)$ be a fixed basis of the lattice $\Gamma$ and denote $e^{2i\pi\theta_j}:=\beta(\gamma_j)$ for $\theta_j\in[0,1[$ and $1\leq j\leq n$.
For all $1\leq j\leq n$ let $\tilde{\beta}(\gamma_j):=e^{i\pi\theta_j}$ and define $\delta_j\in\{0,1\}$ via $(-1)^{\delta_j}:=\varepsilon(\gamma_j)$.
Denoting by $\Gamma^*$ the dual lattice (those $1$-forms on $\RR^n$ with integral values on $\Gamma$) and by $(\sigma_1,\ldots,\sigma_N)$ an orthonormal basis of $\Sigma_n=\CC^N$ (where $N:=2^{\left[\frac{n}{2}\right]}$), we let, for any $\gamma^*\in\Gamma^*$,
\[\theta_{\gamma^*}:=\gamma^*+\frac{1}{2}\sum_{j=1}^n(\delta_j+\theta_j)\gamma_j^*\in(\RR^n)^*,\]
where $(\gamma_1^*,\ldots,\gamma_n^*)$ denotes the dual basis to $(\gamma_1,\ldots,\gamma_n)$.
Fixing a real valued $1$-form $A$ on $\mathbb{T}^n$ or, equivalently, a $\Gamma$-invariant $1$-form on $\RR^n$, we consider the induced connection $1$-form $iA$ on the auxiliary $\mathbb{U}_1$-bundle.
Then the Levi-Civita covariant derivative $\nabla^A$ on spinors induced by the metric and $A$ is given by $\nabla_X^A\phi=\nabla_X\phi+\frac{i}{2}A(X)\phi$, for all vector fields $X$ and spinor fields $\phi$ on $\mathbb{T}^n$.
Here, $\nabla$ denotes the Levi-Civita covariant derivative on spinors, which exists because of the existence of a reduction of the spin$^c$ structure to $\mathrm{Spin}_n\times\mathbb{U}_1$, as explained above.\\

Now let, for $\gamma^*\in\Gamma^*$ and $\ell\in\{1,\ldots,N\}$,
\[\phi_{\gamma^*,\ell}:=e^{2i\pi\theta_{\gamma^*}}\sigma_\ell.\]
This defines a spinor field on $\mathbb{T}^n$: for any $\gamma'\in\Gamma$ and any $x\in\RR^n$,
\begin{eqnarray*}
\phi_{\gamma^*,\ell}(x+\gamma')&=&e^{2i\pi\theta_{\gamma^*}(x+\gamma')}\sigma_\ell\\
&=&e^{2i\pi\theta_{\gamma^*}(\gamma')}\phi_{\gamma^*,\ell}(x)\\
&=&e^{2i\pi(\gamma^*(\gamma')+\frac{1}{2}\sum_{j=1}^n(\delta_j+\theta_j)\gamma_j^*(\gamma'))}\phi_{\gamma^*,\ell}(x)\\
&=&\underbrace{e^{2i\pi\gamma^*(\gamma')}}_{1}\cdot(-1)^{\sum_{j=1}^n\delta_j\gamma_j^*(\gamma')}\cdot e^{i\pi\sum_{j=1}^n\theta_j\gamma_j^*(\gamma')}\phi_{\gamma^*,\ell}(x)\\
&=&\varepsilon(\gamma')\tilde{\beta}(\gamma')\phi_{\gamma^*,\ell}(x).
\end{eqnarray*}
Letting $\theta_{\gamma^*}':=\theta_{\gamma^*}+\frac{A}{4\pi}\in(\RR^n)^*$, we compute, for any vector field $X$ on $\mathbb{T}^n$,
\[\nabla_X^A\phi_{\gamma^*,\ell}=2i\pi\left(\theta_{\gamma^*}(X)+\frac{A(X)}{4\pi}\right)\phi_{\gamma^*,\ell}=2i\pi\theta_{\gamma^*}'(X)\phi_{\gamma^*,\ell}.\]
As a consequence, for the associated Dirac operator $D^A=\sum_{j=1}^n e_j\cdot\nabla_{e_j}^A$, we have
\begin{equation}\label{eq:DAphigammal}
D^A\phi_{\gamma^*,\ell}=2i\pi\theta_{\gamma^*}'\cdot\phi_{\gamma^*,\ell}.
\end{equation}

If $\theta_{\gamma^*}'=0$, which happens iff $A=-4\pi\theta_{\gamma^*}$ (and is in particular constant i.e., parallel), then $\phi_{\gamma^*,\ell}\in\ker(D^A)$.
But in that case $\ker(D^A)$ consists of $\nabla^A$-parallel spinor fields because of the Schr\"odinger-Lichnerowicz formula
\[(D^A)^2=(\nabla^A)^*\nabla^A+\frac{S}{4}\mathrm{Id}+\frac{i}{2}dA\cdot\]
and the fact that $S=0$ (the metric is flat, thus scalar-flat) and $dA=0$ since $A$ is then parallel.
Therefore, if $\theta_{\gamma^*}'=0$, then $\ker(D^A)$ is exactly $N=2^{\left[\frac{n}{2}\right]}$-dimensional and spanned by the $\phi_{\gamma^*,\ell}$, $1\leq \ell \leq N$ (recall that $\gamma^*$ is fixed here). If $\theta_{\gamma^*}'\neq0$, we have to pay attention to the fact that $\theta_{\gamma^*}'$ may vanish pointwise because $A$ is not assumed to be constant.
At those points where $\theta_{\gamma^*}'$ does not vanish, we may split
\[\Sigma_n=\ker\left(i\frac{\theta_{\gamma^*}'}{|\theta_{\gamma^*}'|}\cdot-\mathrm{Id}\right)\bigoplus\ker\left(i\frac{\theta_{\gamma^*}'}{|\theta_{\gamma^*}'|}\cdot+\mathrm{Id}\right).\]
In dimension $n=1$, only one of those subspaces is nonzero, and then complex $1$-dimensional, which forces the other one to vanish.
This means that, still at those points where $\theta_{\gamma^*}'$ does not vanish, $\phi_{\gamma^*,\ell}$ is a $D^A$-eigenspinor associated to the eigenvalue $2\pi|\theta_{\gamma^*}'|$ or $-2\pi|\theta_{\gamma^*}'|$.
In dimension $n\geq2$, both subspaces have to be exactly $\frac{N}{2}$-dimensional since the Clifford action of any (pointwise) nonzero tangent vector which is orthogonal to $\theta_{\gamma^*}'$ (such a vector exists if $n\geq2$) anti-commutes with the Clifford action of $\theta_{\gamma^*}'$ and hence exchanges both subspaces isomorphically.
Then the original constant basis $(\sigma_1,\ldots,\sigma_N)$ of $\Sigma_n$ can be replaced by a pointwise basis $(\sigma_1^+,\ldots,\sigma_{\frac{N}{2}}^+,\sigma_1^-,\ldots,\sigma_{\frac{N}{2}}^-)$ consisting of $\pm1$-eigenvectors of the Clifford action $i\frac{\theta_{\gamma^*}'}{|\theta_{\gamma^*}'|}\cdot$.
The problem is now that the identities above for $\nabla^A$ and hence $D^A$ do not hold any longer since $\sigma_1^+,\ldots,\sigma_{\frac{N}{2}}^+,\sigma_1^-,\ldots,\sigma_{\frac{N}{2}}^-$ now do depend on the base-point and hence are not {\sl a priori} constant, except for instance if $A$ is.
We conclude that this splitting does not allow us to make further progress in the general case which is different to the spin case.\\

If $A$ is parallel, then so is $\theta_{\gamma^*}'$ and, letting $\phi_{\gamma^*,\ell}^\epsilon:=e^{2i\pi\theta_{\gamma^*}}\sigma_{\ell}^\epsilon$ for $1\leq \ell\leq \frac{N}{2}$ and $\epsilon\in\{\pm1\}$, we deduce from \eqref{eq:DAphigammal} that 
\[D^A\phi_{\gamma^*,\ell}^\epsilon=2\pi\epsilon|\theta_{\gamma^*}'|\phi_{\gamma^*,\ell}^\epsilon\]
for all $1\leq \ell\leq \frac{N}{2}$ and both $\epsilon\in\{\pm1\}$.
This provides the spectrum of $D^A$ in that case but beware that it is not necessarily symmetric about $0$ any longer since $-\gamma^*-\sum_{j=1}^n(\delta_j+\theta_j)\gamma_j^*$ does not belong to $\Gamma^*$ any more because of $\theta_j\in[0,1[$ (except when $\theta_j=0$ for all $j$ of course).
\end{proof}

In case $A$ is closed, the Hodge decomposition yields $A=h+df$ for some real function $f$ and harmonic $1$-form $h$, which is actually parallel by the Bochner formula because of the torus being Ricci-flat.
Then, from Proposition \ref{prop:gaugeonva}, the operators $D^A$ and $D^h$ are unitarily equivalent and therefore have the same spectrum.

\begin{remark}\label{rem:diamagneticineqTn}\rm We notice that the diamagnetic inequality for the magnetic Dirac operator may or may not hold on $\mathbb{T}^n$ according to the choice of $\eta$, at least when the underlying spin structure is nontrivial; recall that $0$ lies in the Dirac spectrum when the spin structure on $\mathbb{T}^n$ is trivial i.e., when $\delta_1=\ldots=\delta_n=0$.
Namely if $\theta_1=\ldots=\theta_n=0$ (the auxiliary $\mathbb{U}_1$-bundle should be trivial) and at least one $\delta_j=1$, then the smallest Dirac-eigenvalue in absolute value is given by $\displaystyle |\lambda_1|=2\pi|\gamma^*+\frac{1}{2}\sum_{j=1}^n\delta_j\gamma_j^*|>0$ for some $\gamma^*\in\Gamma^*$ (we keep the notations of Theorem \ref{t:specflattorusspincAparallel}).
But choosing $\eta=\frac{A}{2}$ to be a positive multiple of $\gamma^*+\frac{1}{2}\sum_{j=1}^n\delta_j\gamma_j^*$, the smallest positive $D^{t\eta}$-eigenvalue will be \emph{strictly larger} than $|\lambda_1|$ for small positive $t$.
Conversely, if $\eta=\frac{A}{2}$ is a negative multiple of $\gamma^*+\frac{1}{2}\sum_{j=1}^n\delta_j\gamma_j^*$, then the smallest positive $D^{t\eta}$-eigenvalue will be \emph{strictly smaller} than $|\lambda_1|$ for small positive $t$.
\end{remark}

\section{Killing magnetic field}\label{s:KillingmagneticVF}
In this section, we consider the particular case when the magnetic field is a Killing vector field of constant norm. This gives rise to local Riemannian submersions with one-dimensional fibers given by the integral curves of the magnetic field. 
We will then estimate a part of the spectrum of the corresponding magnetic Dirac operator in terms of the geometry of those submersions. 
In the following, we review some basic facts on spin Riemannian flows, which can be found in \cite{Carriere:84, GinouxHabibspKTgeom:08, HabibRichardson:18,Reinhart:59, Tondeur:88}.\\  

Let $(M^{n},g,\zeta)$ be a closed oriented Riemannian manifold together with a unit Killing vector field $\zeta$, that is $\mathcal{L}_\zeta g=0$.
In this case,  $\zeta$ defines a Riemannian foliation, called Riemannian flow, whose leaves are given by the integral curves of $\zeta$. 
In the following we denote by $Q=\zeta^\perp$ the normal bundle of the flow. 
Locally, a Riemannian flow is given by a Riemannian submersion whose fibers are the leaves of the foliation and the normal bundle corresponds to  the tangent space of the base manifold. 
It is known that the bundle $Q$ carries a natural covariant derivative $\nabla^Q$ given for all $Y\in \Gamma(Q)$ by
\begin{equation*}
\nabla^Q_X Y=
\begin{cases}
	\pi([\zeta,Y]) & \text{if $X=\zeta$,}\\\\
	\pi(\nabla_X Y) & \text{if $X\in \Gamma(Q)$,}
\end{cases}
\end{equation*}  
where $\nabla$ is the Levi-Civita connection on $M$ and $\pi\colon TM\to Q$ is the orthogonal projection. 
This connection $\nabla^Q$ is compatible with the induced metric $g_{|_Q}$ on $Q$ and has a free torsion given by $T^\nabla(X,Y):=\nabla^Q_X \pi(Y)-\nabla^Q_Y \pi(X)-\pi([X,Y])$ for all $X,Y\in \Gamma(TM)$. 
An easy computation gives the relation between $\nabla$ and $\nabla^Q$ through the formulas 
\begin{equation*}
\begin{cases}
	\nabla_\zeta Y=&\nabla^Q_\zeta Y+h(Y),\\\\
	\nabla_X Y=&\nabla^Q_X Y-g(h(X),Y)\zeta,
\end{cases}
\end{equation*}  
for all $X,Y\in \Gamma(Q)$, where $h:=\nabla\zeta$ is a skew-symmetric endomorphism on $TM$ called the O'Neill tensor \cite{ONeill:66}. 
Now assume that $M$ is spin and let us denote by $\Sigma M$ its spinor bundle.
Since $TM=\RR\zeta\oplus Q$, the bundle $Q$ carries also a spin structure and its spinor bundle $\Sigma Q$ can be canonically identified with the one on $M$ when $n$ is odd and when $n$ is even, we have $\Sigma M\simeq \Sigma Q\oplus \Sigma Q$. 
Also the Clifford multiplication on $M$ and $Q$ can be identified by $Y\cdot_Q\varphi=Y\cdot\varphi$ when $n$ is odd and, when $n$ is even, we have $Y\cdot\zeta\cdot\varphi=(Y\cdot_Q\oplus-Y\cdot_Q)\varphi$ for all $Y\in \Gamma(Q)$. 
The Clifford action of $i\zeta$ is given by ${\rm Id}_{\Sigma^+ Q}\oplus -{\rm Id}_{\Sigma^-Q}$ if $n$ is odd and by $i\zeta\cdot=\left(\begin{array}{cc}0&\mathrm{Id}_{\Sigma Q}\\\mathrm{Id}_{\Sigma Q}&0\end{array}\right)$ if $n$ is even. 
The spinorial connections of $\Sigma Q$ and $\Sigma M$ are related by 
\begin{equation}\label{eq:oneilforspin}
\begin{cases}
	\nabla_\zeta\varphi=&\nabla^{\Sigma Q}_\zeta \varphi+\frac{1}{2}\Omega\cdot \varphi, \\\\
	\nabla_Y\varphi=&\nabla^{\Sigma Q}_Y \varphi+\frac{1}{2}\zeta\cdot h(Y)\cdot\varphi,
\end{cases}
 \end{equation}
for all $Y\in \Gamma(Q)$. 
Here $\Omega:=\frac{1}{2}d\zeta^\flat$ is the two-form that is associated to $h$ via $\Omega(Y,Z):=g(h(Y),Z)$, for any $Y,Z\in \Gamma(Q)$. 
Notice that the notion of the Lie derivative of any spinor in the direction of $\zeta$ defined in \cite{BourguignonGauduchon:92} is just the covariant derivative $\nabla^{\Sigma Q}_\zeta$. 
Namely, by \cite[Prop. 17]{BourguignonGauduchon:92}, the Lie derivative of a spinor field $\varphi$ is expressed by the formula
$\nabla_\zeta\varphi=\mathcal{L}_\zeta\varphi+\frac{1}{4}d\zeta^\flat\cdot\varphi$. 
Hence using that $d\zeta^\flat=2\Omega$ and by comparing with \eqref{eq:oneilforspin}, we get $\mathcal{L}_\zeta\varphi=\nabla^{\Sigma Q}_\zeta\varphi$. 
Therefore basic spinors, i.e. those spinors $\varphi$ satisfying $\nabla^{\Sigma Q}_\zeta\varphi=0$, correspond to the so-called projectable spinors. 
The transversal Dirac operator is the first-order differential operator defined by $D_Q:=\sum_{k=1}^{n-1} e_k\cdot_Q\nabla^{\Sigma Q}_{e_k}$ on $\Gamma(\Sigma Q)$, where $\{e_k\}_{k=1,\ldots,n-1}$ is a local orthonormal frame of $\Gamma(Q)$. 
It is a transversally elliptic and self-adjoint operator when restricted to basic spinors, in particular it has a discrete spectrum \cite{ElKacimi:90}. 
This is called the basic Dirac operator and is often denoted by $D_b$. 
Now, with the help of \eqref{eq:oneilforspin}, the Dirac operator $D$ on $M$ is related to the transversal Dirac operator $D_Q$ by \cite{HabibRichardson:18} 
\begin{equation}\label{eq:diracrelation}
D=\left\{
\begin{matrix}
D_Q-\frac{1}{2}\zeta\cdot\Omega\cdot+\zeta\cdot\nabla^{\Sigma Q}_\zeta & \text{if $n$ is odd,}\\\\
	\zeta\cdot(D_Q\oplus (-D_Q))-\frac{1}{2}\zeta\cdot\Omega\cdot+\zeta\cdot\nabla^{\Sigma Q\oplus\Sigma Q}_\zeta & \text{if $n$ is even.}
\end{matrix}\right.
\end{equation}  
It is shown in \cite[Lem. 2.6]{HabibRichardson:18} that the Dirac operator $D$ preserves the set of basic spinors $\Gamma_b(\Sigma Q)$ and, thus, it decomposes as an $L^2$-orthogonal sum $D_{|_{\Gamma_b(\Sigma Q)}}\oplus D_{|_{\Gamma_b(\Sigma Q)^\perp}}$ when $n$ is odd. 
An analogous decomposition holds when $n$ is even. 
Therefore, the spectrum of $D$ consists of eigenvalues of the form $\{\lambda_j\}_{j=1}^\infty\cup \{\mu_k\}_{k=1}^\infty$ that correspond to the restriction of $D$ to $L^2(\Gamma_b(\Sigma Q))$ and $L^2(\Gamma_b(\Sigma Q))^\perp$. \\

We set $\eta=\zeta^\flat$ to be the one-form on $M$ associated to the vector field $\zeta$ by the musical isomorphism and consider the magnetic Dirac operator $D^{t\eta}=D+it\eta\cdot$, for $t\in \mathbb{R}$. 
Since by \cite[p.71]{GinouxHabibspKTgeom:08} (see also \cite[Lem. 2.2]{HabibRichardson:18}) we have $[\nabla^{\Sigma Q}_X,\eta\cdot]=0$, for any vector field $X\in \Gamma(TM)$, the magnetic Dirac operator preserves the set of basic spinors as well as its $L^2$-orthogonal complement. 
Therefore, we shall denote by $\{\lambda^{t\eta}_j\}_{j=1}^\infty\cup \{\mu^{t\eta}_k\}_{k=1}^\infty$ the set of eigenvalues corresponding to this decomposition. 
In the following, we will give estimates for the eigenvalues $\{\lambda^{t\eta}_j\}_{j=1}^\infty$. 
For this, we need the following lemma:

\begin{lemma}\label{l:muomegavarphi} Let $(M^{n},g,\zeta)$ be a closed Riemannian spin manifold of odd dimension $n$ equipped with a unit Killing vector field $\zeta$. 
Let $\varphi$ be a basic eigenspinor of the magnetic Dirac operator $D^{t\eta}$ associated with an eigenvalue $\lambda^{t\eta}$, we have 
\begin{equation}\label{eq:muomegavarphi}
\lambda^{t\eta}\int_M \Re\langle i\zeta\cdot\varphi,\varphi\rangle d\mu_g=-\frac{1}{2}\int_M\langle\Omega\cdot\varphi,i\varphi\rangle d\mu_g+t\int_M|\varphi|^2 d\mu_g.
\end{equation}
\end{lemma}
\begin{proof} 
Using the fact that $[\nabla^{\Sigma Q}_X,\eta\cdot]=0$ for any $X\in \Gamma(TM)$ and the identification of the Clifford multiplications between $M$ and $Q$ as $n$ is odd, we can easily deduce that $D_b(\zeta\cdot)=-\zeta\cdot D_b$. 
We compute for any basic spinor $\varphi$
\begin{eqnarray*}
    \int_M \langle D_b\varphi,i\zeta\cdot\varphi\rangle d\mu_g&=&\int_M  \langle \varphi,iD_b(\zeta\cdot\varphi)\rangle d\mu_g\\
    &=&-\int_M  \langle \varphi,i\zeta\cdot D_b\varphi\rangle d\mu_g\\
    &=&-\int_M  \langle i\zeta\cdot\varphi,D_b\varphi\rangle d\mu_g.
\end{eqnarray*}
Hence, we deduce that $\int_M\Re\langle D_b\varphi,i\zeta\cdot\varphi\rangle d\mu_g=0$ for any spinor field $\varphi$. When $\varphi$ is an eigenspinor of the magnetic Dirac operator, we use \eqref{eq:diracrelation} to deduce \eqref{eq:muomegavarphi}.
\end{proof}
Now, we state the main result of this section:

\begin{theorem}\label{thm:estimatebasic}
Let $(M^{n}, g,\zeta)$ be a closed Riemannian spin manifold with nonnegative scalar curvature $S$. 
We assume that $M$ carries a unit Killing vector field $\zeta$. 
Then any eigenvalue of the Dirac operator $D^{t\eta}$ restricted to basic spinors (or equivalently projectable spinors) satisfies, for $n>3$,
$$|\lambda^{t\eta}|\geq \mathop{\rm inf}\limits_M\left(-\frac{[\frac{n-1}{2}]^\frac{1}{2}|\Omega|}{2}+\sqrt{t^2+\frac{n-1}{4(n-2)}(S+2|\Omega|^2)}\right)$$ 
and, for $n=3$, the first positive eigenvalue $\lambda_1^{t\eta}$ satisfies  
$$\lambda_1^{t\eta}\geq \mathop{\rm inf}\limits_M\left(\frac{b}{2}+\sqrt{t^2+\frac{1}{2}(S+2b^2})\right),$$
where $h=bJ$ is the O'Neill tensor.
\end{theorem}
\begin{proof} The key point of the proof is to define the transversal twistor operator as follows: For all $X\in \Gamma(Q)$ and $\varphi\in \Gamma_b(\Sigma Q)$ 
$$P_X^Q\varphi:=\nabla^{\Sigma Q}_X\varphi+\frac{1}{n-1}X\cdot_Q D_b\varphi.$$
An easy computation shows that $|P^Q\varphi|^2=|\nabla^{\Sigma Q}\varphi|^2-\frac{1}{n-1}|D_b\varphi|^2.$
Therefore, with the help of the transversal Schr\"odinger-Lichnerowicz identity $D_b^2=\nabla^*\nabla+\frac{1}{4}S^Q$ where $S^Q$ is the transversal scalar curvature \cite{GlazebrookKamber:91}, we deduce that 
\begin{equation}\label{eq:transtwistor}
\int_M|P^Q\varphi|^2 d\mu_g=\frac{n-2}{n-1}\int_M|D_b\varphi|^2 d\mu_g-\frac{1}{4}\int_M S^Q|\varphi|^2 d\mu_g.
\end{equation}
Let us first consider the case when $n$ is odd. 
Let $\varphi$ be a basic spinor which is also eigenspinor for $D^{t\eta}$ associated with the eigenvalue $\lambda^{t\eta}$. 
Identity \eqref{eq:transtwistor} reduces to 
\begin{align*}
    0&\leq  \frac{n-1}{n-2}\int_M|P^Q\varphi|^2 d\mu_g \\
    &=\int_M|D_b\varphi|^2 d\mu_g-\frac{n-1}{4(n-2)}\int_M S^Q |\varphi|^2 d\mu_g\\
    &\stackrel{\eqref{eq:diracrelation}}{=}\int_M\bigg((\lambda^{t\eta})^2|\varphi|^2+\frac{1}{4}|\Omega\cdot\varphi|^2+t^2|\varphi|^2-\lambda^{t\eta}\Re(\langle\Omega\cdot\varphi,\zeta\cdot\varphi\rangle)-2t\lambda^{t\eta}\Re (i\langle\zeta\cdot\varphi, \varphi\rangle)\\&\hspace{1cm}-t\Re\langle\Omega\cdot\varphi,i\varphi\rangle\bigg) d\mu_g-\frac{n-1}{4(n-2)}\int_M S^Q |\varphi|^2d\mu_g\\
    &\stackrel{\eqref{eq:muomegavarphi}}{=} \int_M\left((\lambda^{t\eta})^2|\varphi|^2+\frac{1}{4}|\Omega\cdot\varphi|^2-t^2|\varphi|^2-\lambda^{t\eta}\Re(\langle\Omega\cdot\varphi,\zeta\cdot\varphi\rangle)\right)d\mu_g-\frac{n-1}{4(n-2)}\int_M S^Q |\varphi|^2d\mu_g.
   \end{align*}
The case when $n$ is even gives also the same inequality as above. 
Namely, the spinor bundle of $M$ can be identified with $\Sigma Q\oplus \Sigma Q$ and we write any eigenspinor of $D^{t\eta}$ as $\varphi=\varphi_1+\varphi_2$. 
Now the equation $D^{t\eta}\varphi=\lambda^{t\eta}\varphi$ can be equivalently written as 
$$D\varphi_1=\lambda^{t\eta}\varphi_2-it\zeta\cdot\varphi_1\quad\text{and}\quad D\varphi_2=\lambda^{t\eta}\varphi_1-it\zeta\cdot\varphi_2.$$
Therefore, using \eqref{eq:diracrelation}, we get the following expressions for the basic Dirac operator of $\varphi_1$ and $\varphi_2$ 
\begin{equation}\label{eq:basiceven}
D_b\varphi_1=-\lambda^{t\eta}\zeta\cdot\varphi_2+\frac{1}{2}\Omega\cdot\varphi_1-it\varphi_1\quad\text{and}\quad D_b\varphi_2=\lambda^{t\eta}\zeta\cdot\varphi_1-\frac{1}{2}\Omega\cdot\varphi_2+it\varphi_2.
\end{equation}
Applying Inequality \eqref{eq:transtwistor} to $\varphi_1$ gives that 
\begin{align}
\label{eq:inequalityeven}
0\leq &\int_M\bigg((\lambda^{t\eta})^2|\varphi_2|^2+\frac{1}{4}|\Omega\cdot\varphi_1|^2+t^2|\varphi_1|^2-\lambda^{t\eta}\Re(\langle\Omega\cdot\varphi_1,\zeta\cdot\varphi_2\rangle) \\
&\hspace{1cm}+2t\lambda^{t\eta}\Im (\langle\zeta\cdot\varphi_2,\varphi_1\rangle)-t\Im(\langle\Omega\cdot\varphi_1,\varphi_1\rangle)\bigg) d\mu_g\nonumber-\frac{n-1}{4(n-2)}\int_M S^Q |\varphi_1|^2d\mu_g.
\end{align}
By taking the Hermitian inner product of the first equation in \eqref{eq:basiceven} with $\varphi_1$ and identifying the imaginary parts, we get that 
$$\lambda^{t\eta}\int_M\Im(\langle\zeta\cdot\varphi_2,\varphi_1\rangle) d\mu_g-\frac{1}{2}\int_M\Im  (\langle \Omega\cdot\varphi_1,\varphi_1\rangle) d\mu_g=-t\int_M|\varphi_1|^2 d\mu_g.$$
Replacing this last identity into \eqref{eq:inequalityeven} gives that 
\begin{align*}
0&\leq \int_M\left((\lambda^{t\eta})^2|\varphi_2|^2+\frac{1}{4}|\Omega\cdot\varphi_1|^2-t^2|\varphi_1|^2-\lambda^{t\eta}\Re(\langle\Omega\cdot\varphi_1,\zeta\cdot\varphi_2\rangle)\right)-\frac{n-1}{4(n-2)}\int_M S^Q |\varphi_1|^2d\mu_g.
\end{align*}
Now, we do the same computations as before by applying inequality \eqref{eq:transtwistor} to $\varphi_2$ and find after adding both inequalities 
\begin{align*}
0\leq \int_M\left((\lambda^{t\eta})^2|\varphi|^2+\frac{1}{4}|\Omega\cdot\varphi|^2-t^2|\varphi|^2-\lambda^{t\eta}\Re(\langle\Omega\cdot\varphi,\zeta\cdot\varphi\rangle)\right)d\mu_g-\frac{n-1}{4(n-2)}\int_M S^Q |\varphi|^2d\mu_g.
\end{align*}
  
In the following, we will distinguish the cases when $n=3$ and $n>3$. 
When $n=3$, the O'Neill tensor $h$ can be written as $h=bJ$ for some real valued function $b:M\to \RR$, where $J$ is the complex structure on the bundle $Q$. 
In this case, we have $\Omega\cdot\varphi=b\zeta\cdot\varphi$.
Hence the above inequality reduces to 
\begin{align*}
0&\leq \int_M\left((\lambda^{t\eta})^2+\frac{b^2}{4}-t^2-b\lambda^{t\eta}-\frac{1}{2} S^Q \right)|\varphi|^2d\mu_g\\
&=\int_M \left(\lambda^{t\eta}-\frac{b}{2}+\sqrt{t^2+\frac{1}{2}S^Q}\right)\left(\lambda^{t\eta}-\frac{b}{2}-\sqrt{t^2+\frac{1}{2}S^Q}\right)|\varphi|^2 d\mu_g.
\end{align*}
Recall now that the following relation $S^Q=S+2|\Omega|^2$ holds \cite{ONeill:66}. 
Thus, if $S\geq 0$, then $S^Q$ is also nonnegative and we can write 
$$\sqrt{t^2+\frac{1}{2}S^Q}=\sqrt{t^2+\frac{1}{2}S+|\Omega|^2}\geq |\Omega|=|b|\geq \pm \frac{b}{2}.$$
Hence the first positive eigenvalue satisfies 
$$\lambda_1^{t\eta}\geq \mathop{\rm inf}\limits_M\left(\frac{b}{2}+\sqrt{t^2+\frac{1}{2}S^Q}\right).$$
When $n>3$, we use the pointwise estimate $|\Omega\cdot\varphi|\leq [\frac{n-1}{2}]^\frac{1}{2}|\Omega||\varphi|$ which can be proven in the same way as in \cite[Lem. 3.3]{HerzlichMoroianu:99} to obtain
\begin{align*}
0&\leq \int_M\left((\lambda^{t\eta})^2+\frac{1}{4}[\frac{n-1}{2}]|\Omega|^2-t^2+|\lambda^{t\eta}|[\frac{n-1}{2}]^\frac{1}{2}|\Omega|-\frac{n-1}{4(n-2)} S^Q \right)|\varphi|^2d\mu_g\\
&= \int_M \left(|\lambda^{t\eta}|+\frac{[\frac{n-1}{2}]^\frac{1}{2}|\Omega|}{2}+\sqrt{t^2+\frac{n-1}{4(n-2)}S^Q}\right)\left(|\lambda^{t\eta}|+\frac{[\frac{n-1}{2}]^\frac{1}{2}|\Omega|}{2}-\sqrt{t^2+\frac{n-1}{4(n-2)}S^Q}\right)|\varphi|^2 d\mu_g,
\end{align*}
which gives the estimate 
$$|\lambda^{t\eta}|\geq \mathop{\rm inf}\limits_M\left(-\frac{[\frac{n-1}{2}]^\frac{1}{2}|\Omega|}{2}+\sqrt{t^2+\frac{n-1}{4(n-2)}S^Q}\right)$$
as required. 
Notice that, when the parameter $t$ is big enough, this lower bound is positive. 
\end{proof}

\begin{remark}\label{rem:estimatebasicS3}\rm When $M=\mathbb{S}^3$ equipped with the metric of curvature $1$, we take $b=1$ and, therefore, the estimate in Theorem \ref{thm:estimatebasic} becomes 
$$\lambda_1^{t\eta}\geq \frac{1}{2}+\sqrt{t^2+4}.$$
This lower bound appears in the set of eigenvalues studied in Theorem \ref{t:spectrumDetaS3}. 
Also, on the round sphere $(\mathbb{S}^{n},g,\eta)$ with $n=2m+1>3$ and $\eta$ is the Reeb vector field, the lower bound obtained in \eqref{eq:Hijazimagneticineq} is  $|\lambda^{t\eta}|\geq\frac{n}{2}-t$. 
However, the lower bound in Theorem \ref{thm:estimatebasic} is 
$$|\lambda^{t\eta}|\geq\left(-\frac{n-1}{4}+\sqrt{t^2+\frac{(n-1)^2(n+1)}{4(n-2)}}\right).$$
Hence, for $t\geq \frac{3n-1}{4}$ the above lower bound is better than $\frac{n}{2}-t$. 
\end{remark}

\end{document}